\def\versiondate{5 July 2005}
\input math.macros
\input Ref.macros

%\proofmodetrue
%\leftsectionheadtrue
\checkdefinedreferencetrue
%\continuousnumberingtrue
\continuousfigurenumberingtrue
\theoremcountingtrue
\sectionnumberstrue
%\figuresectionnumberstrue
\forwardreferencetrue
%\lefteqnumberstrue
%\tocgenerationtrue
\citationgenerationtrue
\nobracketcittrue
\hyperstrue
\initialeqmacro

\input\jobname.key
\bibsty{myapalike}

\def\cp{\tau}   % complexity
\def\gh{G}  % graph
\def\gp{\Gamma}  % group
\def\gpe{\gamma}  % group element
\def\verts{{\ss V}}
\def\vertex{{\ss V}}
\def\edges{{\ss E}}
\def\tr{{\rm tr}}  % usual trace
\def\Tr{{\rm Tr}}  % normalized trace
\def\Det{{\rm Det}} % normalized det
\def\bfo{{\bf 1}}
\def\bfz{{\bf 0}}
\def\detp{{\rm det}'}
\def\fo{{\frak F}}  %% forest
  %% tree
\def\wsf{{\ss WSF}}
\def\fsf{{\ss FSF}}

\def\bdv{\partial_{\verts}}
\def\cpe{\cp_{\rm e}}
\def\Aut{{\rm Aut}} % automorphism group
\def\grn{{\cal G}}  % Green gen function = return series
\def\rad{{\bf r}}   % radius of cnv
\def\cat{{\rm G}}  % Catalan's constant
\def\asym{{\bf h}}   % asymptotics
  % local entropy
\def\ip#1{(\changecomma #1)_{\pi}}
\def\bigip#1{\bigl(\changecomma #1\bigr)_\pi}

\def\bigpip#1{\bigl(\changecomma #1\bigr)}
\def\Bigpip#1{\Bigl(\changecomma #1\Bigr)}
\def\changecomma#1,{#1,\,}
\def\bigchangecomma#1,{#1,\;}
\def\leftchangecomma#1,{#1,\ }
\def\G{{\cal G}}
\def\PGW{{\ss PGW}}
\def\sgn{{\rm sgn}}
\def\rtd{\rho}  % the measure on rooted graphs
\def\bp{o}
\def\one{{\bf 1}}
\def\bfz{{\bf 0}}
\def\T{[0, 1]}
\def\ent{{\bf H}}  % entropy
\def\alg{{\ss Alg}}  % algebra
\def\expdeg{\overline{\rm deg}}  % expected degree

\def\BLPSusf{\ref b.BLPS:usf/, hereinafter referred to as BLPS (2001)%
\def\BLPSusf{BLPS \htmllocref{\bibcode{BLPS:usf}}{(2001)}}}

\def\BLPSgip{\ref b.BLPS:gip/, hereinafter referred to as BLPS (1999)%
\def\BLPSgip{BLPS \htmllocref{\bibcode{BLPS:gip}}{(1999)}}}

\ifproofmode \relax \else\head{To appear in {\it Combin. Probab. Comput.},
\copyright Cambridge University Press}
{Version of \versiondate}\fi 
\vglue20pt

\title{Asymptotic Enumeration of Spanning Trees}

\author{Russell Lyons}

\abstract{We give new general formulas for the
asymptotics of the number of spanning trees of a large graph.
A special case answers a question of McKay (1983) for regular graphs.
The general answer involves a quantity for infinite graphs that we call
``tree entropy", which we show is a logarithm of a normalized
determinant of the graph Laplacian for infinite graphs.
Tree entropy is also expressed using random walks.
We relate tree entropy to the metric
entropy of the uniform spanning forest process on
quasi-transitive amenable graphs,
extending a result of Burton and Pemantle (1993).
}

\bottomII{Primary 
 05C05, %   Trees
60C05. % Combinatorial probability
Secondary
60K99,  % special processes
 05C80, %   Random graphs
82B20, %lattice systems (Ising, dimer, Potts, etc.) and systems on graphs
28D05, %measure-preserving transformations
37A05. %measure-preserving transformations
}
{Asymptotic complexity, graphs, domino tilings, tree entropy, random walks,
metric entropy, uniform spanning forests, determinant, Laplacian, trace.}
{Research partially supported by NSF grants DMS-0103897 and DMS-0231224.}

\bsection{Introduction}{s.intro}

Methods of enumeration of spanning trees in a finite graph $\gh$ and relations
to various areas of mathematics and physics have been investigated for more
than 150 years.
The number of spanning trees is often called the {\bf complexity} of the
graph, denoted here by $\cp(\gh)$.
The best known formula for the complexity, proved in every basic text on
graph theory, is called ``the Matrix-Tree Theorem", which expresses it as a
determinant.
One is often interested in asymptotics of the complexity of a sequence of
finite graphs that approach (in some sense) some infinite graph.
Use of the Matrix-Tree
Theorem typically involves calculating eigenvalues and their
asymptotics.\Fnote{It appears not to have been noticed before that
in the case of Euclidean lattice graphs, one can dispense with
eigenvalues and go directly to the limit by invoking 
Szeg\H{o}'s limit theorem. In one dimension, this theorem appears in, e.g.,
\ref b.GZ/, p.~44. In higher dimensions, it is due to \ref b.HelLow/, but the
form in which they present it does not involve determinants. See the last
part of the proof of Theorem 5.11 in \ref b.LyonsSteif:dyn/ to see how to
transform their result to give the asymptotics of determinants. In any case,
use of our formula is simpler still.}
This often leads to a formula of the form 
$$
\lim_{n \to\infty} {1 \over |\verts(\gh_n)|}  \log \cp(\gh_n) 
=
\int \log f
$$
for some function $f$ on the unit cube of some Euclidean
space, where $\verts(\gh_n)$ is the vertex set of $\gh_n$.
One of the cases that is most well known, due to its connection with domino
tilings, is that where the graphs $\gh_n$ approach the usual square lattice,
i.e., the nearest-neighbor graph on $\Z^2$. In this case, one has that
$$
\lim_{n \to\infty} |\verts(\gh_n)|^{-1} \log \cp(\gh_n) 
=
\int_0^1 \int_0^1 \log (4 - 2 \cos 2\pi x - 2 \cos 2\pi y) \,dx \,dy
=
4\cat/\pi
\,,
$$
where $\cat := \sum_{k=0}^\infty (-1)^k / (2k+1)^2 = 0.9160^- $
is {\bf Catalan's constant}.
See, e.g., \ref b.BurPem/ and \ref b.ShrockWu/ and the references therein
for this and several other such examples.
We shall reformulate the Matrix-Tree Theorem as an infinite series whose
terms have a probabilistic meaning (see \ref p.finite/).
In fact, the terms of the series are ``local" quantities.
This will permit us to find the asymptotic complexity in a very general
setting (see \ref t.asym/). % and the end of \ref s.asym/).
Other than assuming a bound on the average degree of the finite graphs
of whose complexity we find the asymptotics,
our results are essentially the most general possible. 
In the simplest case, where the finite graphs $\gh_n$ tend to a fixed
transitive infinite graph $\gh$, we shall prove that 
$$
\lim_{n \to\infty} |\verts(\gh_n)|^{-1} \log \cp(\gh_n)
=
\asym(\gh)
:=
\log \deg_\gh(o) - \sum_{k \ge 1} p_k(o;\gh)/k
\,,
$$
where $o$ is a fixed vertex of $\gh$ and $p_k(o;\gh)$ is the probability that
simple random walk started at $o$ on $\gh$ is again at $o$ after $k$ steps.
We term $\asym(\gh)$ the {\bf tree entropy}\Fnote{Another possible
name, ``combinatorial entropy", is already in use with
a variety of meanings.} of $\gh$.
We have chosen this terminology to reflect both its provenance as a limit of
entropies per vertex of uniform spanning trees of finite graphs as well as
its agreement (up to normalization) with the metric entropy of uniform
spanning forests in certain infinite graphs (see \ref t.amen-ent/), as we
shall discuss shortly.

%This asymptotic is also expressed as an infinite series whose terms
%involve probabilities on the limiting graph, rather
%than an integral as above.
%Alternatively, one could express our series as an integral over $[0,
%1]$ of a modification of the return probability generating function, as in \ref
%e.intgrn/ below. This is different from other integral
%representations in the literature.
%This has several advantages: In many cases, there is no known natural
%integral representation, nor any reason to believe that one
%exists. 
%Indeed, the integrals arise only when $\Z^d$ acts quasi-transitively on the
%limiting graph.

Our result allows us to
find asymptotics that were not known previously.
For a simple example,
in case the graphs $\gh_n$ approach the regular tree of degree 4,
then 
$$
\lim_{n \to\infty} |\verts(\gh_n)|^{-1}  \log \cp(\gh_n) 
=
3 \log (3/2)
\,;
$$
see \ref x.wired-free/.
This answers a question of \ref b.McKay/, who had shown such a result only
under stronger conditions on the sequence of graphs $\gh_n$, and who had
noted (in the case of regular graphs) that, if sufficient, our conditions
would be weakest possible.

In the past, it has been observed by direct calculation
that one has the same asymptotic
complexity both of rectangular portions of $\Z^d$ as well as of
$d$-dimensional tori.
Our first main result, \ref t.asym/, shows that this stability phenomenon is
quite general; see \ref c.subgraph/.

Our infinite series representation of the asymptotic complexity has an
additional benefit.
Namely, even when an integral representation is known, it turns out that
for numerical estimation, our infinite series can be much better than the
integral\Fnote{In the $\Z^d$-transitive case, it is not hard to pass
directly between our series representation and the standard integral
representation.}; see, e.g., \ref b.FelkerLyons/ for examples comparing the
two approaches.

In order to state our second main result, we recall the Matrix-Tree
Theorem.
Given a graph $\gh$, define the following matrices indexed by
the vertices of $\gh$:
Let $D_\gh$ be the diagonal matrix whose $(x, x)$-entry equals the degree
of $x$; this is the {\bf degree matrix}.
Let $A_\gh$ be the matrix whose $(x, y)$-entry equals the number of edges
joining $x$ and $y$; this is the {\bf adjacency matrix}.
Finally, let $\Delta_\gh := D_\gh - A_\gh$, the {\bf graph Laplacian
matrix}.
The Matrix-Tree Theorem says that given any finite graph $\gh$,
every cofactor of $\Delta_\gh$ 
is the same, namely, $\cp(\gh)$ (see, e.g., \ref b.GodRoy/, Lemma 13.2.3).
Another version of the Matrix-Tree Theorem states that when $\gh$ is
connected,
$$
\cp(\gh) = |\verts(\gh)|^{-1}
\detp \Delta_\gh
\,,
\label e.MTT
$$
where $\detp A$ denotes the product of the
non-zero eigenvalues of a matrix $A$
%(see, e.g., \ref b.ChungYau:cover/ or \ref b.CDS/).
%(to see this, calculate the coefficient of $t$ in $\det (\Delta_\gh - t I)$ in
%two ways, one way by using the eigenvalues of $\Delta_\gh$ and the other by
%using co-factors of the diagonal and the Matrix-Tree Theorem).
(see, e.g., %\ref b.CDS/ \msnote{Proposition 1.3 in 1980 ed.} or
\ref b.GodRoy/, Lemma 13.2.4).
Thus, the asymptotic complexity of a sequence of graphs is a
limit of the logarithm of the determinant of the graph Laplacians,
appropriately normalized.
Our second main aim is to give meaning to and to prove the
formula 
$$
\asym(\gh) = \log \Det \, \Delta_\gh
$$
directly in terms
of a normalized determinant of the graph Laplacian for infinite graphs.
We shall also use this formula to prove inequalities for tree entropy and to
calculate easily and quickly the classical tree entropy for Euclidean
lattices.
Indeed, this result, \ref t.trlog/, provides a quick way to derive
virtually all the prior asymptotics of this type in the literature, while
\ref t.asym/ gives the rest (and more).

%\msnote{Finer asymptotics; cf Duplantier.}

We have alluded to the fact that tree entropy
also arises as the entropy per vertex of a measure that
is the weak limit of the uniform measures on spanning trees of finite graphs.
To state this precisely, we must first recall the work of
\ref b.Pemantle:ust/.
He proved a conjecture of the present author, namely,
that if an infinite connected graph $\gh$ is exhausted by a sequence of finite
connected subgraphs $\gh_n$, then the weak limit of the uniform spanning tree
measures on $\gh_n$ exists. 
However, it may happen that the limit measure is not supported on trees,
but on forests.
This limit measure is now called the {\bf free uniform spanning forest} on
$\gh$, denoted $\fsf$, or $\fsf_\gh$ when we want to indicate the graph
$\gh$. 
If $\gh$ is itself a tree, then this measure is
trivial, namely, it is concentrated on $\{\gh\}$. Therefore, \ref
b.Hag:umesft/
introduced another limit that had been considered on $\Z^d$
more implicitly by \ref b.Pemantle:ust/ and explicitly by
\ref b.Hag:rcust/, namely, the weak limit of the uniform
spanning tree measures on $\gh_n^*$, where $\gh_n^*$ is the graph $\gh_n$ with
its boundary identified to a single vertex. As \ref b.Pemantle:ust/ showed,
this limit also always exists
on any graph and is now called the {\bf wired uniform spanning forest},
denoted $\wsf$ or $\wsf_\gh$.  
It is clear that both $\fsf$ and $\wsf$ are concentrated on the set of
spanning forests\Fnote{By a
``spanning forest'', we mean a subgraph without cycles that contains every
vertex.} of $\gh$ that are {\bf essential}, meaning that all their
trees are infinite.
Furthermore, since the limits exist regardless of the exhaustion chosen, both
$\fsf_\gh$ and $\wsf_\gh$ are invariant under any automorphisms that $\gh$ may
have.
As shown by \ref b.Hag:rcust/, $\fsf_\gh = \wsf_\gh$ on every amenable
transitive graph $\gh$ such as $\Z^d$.
In fact, the proof of this result given by \BLPSusf, is easily modified
to show that the same holds for every quasi-transitive amenable graph, using
the notion of natural frequency distribution recalled below in \ref s.ent/.
Both $\fsf$ and $\wsf$ are important in their own right; see \ref
b.Lyons:bird/ for a survey and \BLPSusf\ for a comprehensive
treatment.

Suppose now that $\gp$ is an amenable group that acts quasi-transitively on a
graph $\gh$.
Since $\fsf_\gh$ is defined as a limit of uniform measures, it is natural to
expect that it has maximum (metric) entropy in an appropriate class of
$\gp$-invariant measures. Since the set of essential spanning forests is
closed, one generally
considers the class of $\gp$-invariant measures on this set.
Furthermore, one might expect that the entropy of $\fsf_\gh$ is related in
a simple way to the exponential growth rate of the complexity of the finite
subgraphs $\gh_n$ exhausting $\gh$.
Finally, one might wonder whether $\fsf_\gh$ is unique as a measure of maximal
entropy on the essential spanning forests.

The case of $\gp = \Z^d$ was considered by \ref b.BurPem/.
Their Theorem 6.1 gave a positive answer to all three questions implicit in
the preceding paragraph.
However, it turns out that one of the claims, the uniqueness of the measure
of maximal entropy, received an incorrect proof.
The correct proofs of the other claims
relied to some extent on the natural tiling of $\Z^d$ by
large boxes, which is not available in general amenable quasi-transitive
graphs.
The third aim of this paper is to extend these two other claims 
to amenable quasi-transitive graphs, that is, to prove that the formula for
entropy is correct and to prove that this entropy is maximal (see \ref
t.amen-ent/).
We have not been able to prove uniqueness of the measure of maximal
entropy.
Therefore, the deduction by \ref b.BurPem/ of the uniqueness of the measure of
maximal entropy for domino tilings in $\Z^2$ (their Theorem 7.2)
remains incomplete.
However, a new and more general proof for the uniqueness of the measure of
maximal entropy for domino tilings has been found by \ref b.Sheff:thesis/.
The arguments in \ref b.BurPem/ can then be used to deduce the uniqueness
of the measure of maximal entropy for spanning trees in $\Z^2$ (but not in
higher dimensions, which remains open).

We give general background and definitions related to graphs, groups, and
entropy in \ref s.back/.
In \ref s.asym/, we give our reformulation of the Matrix-Tree Theorem and its
asymptotics, together with several examples.
We give the formula for tree entropy as a determinant on infinite graphs
in \ref s.further/, with applications.
We prove the assertions on metric entropy in \ref s.ent/.

{\bf Note added in proof:} It has just been discovered that there is an
error in \ref b.AL:uni/. Namely, there is a gap in the proof
that all unimodular random
rooted graphs are random weak limits of finite graphs. This makes
the proofs of certain results here incomplete. 
\ref r.explicit/, \ref p.positive/, and \ref t.zero/ will be
complete when restricted to random weak
limits of finite graphs, rather than claimed for all unimodular random
rooted graphs.

\bsection{Background}{s.back}

We denote a (multi-)graph $\gh$ with vertex set $\verts$ and edge set $\edges$
by $\gh = (\verts, \edges)$.
When there is more than one graph under discussion, we write $\verts(\gh)$ or
$\edges(\gh)$ to avoid ambiguity.
We denote the degree of a vertex $x$ in a graph $\gh$ by $\deg_\gh(x)$.
Unless stated otherwise, we assume all degrees finite.
{\bf Simple random walk} on $\gh$ is the Markov chain whose state space is
$\verts$ and whose transition probability from $x$ to $y$ equals the number of
edges joining $x$ to $y$ divided by $\deg_\gh(x)$.

An infinite path in a tree that starts at any vertex and does not
backtrack is called a {\bf ray}. Two rays are {\bf equivalent} if they have
infinitely many vertices in common. An equivalence class of rays is
called an {\bf end}. 
More generally, an {\bf end} of an infinite graph $\gh$ is an equivalence
class of infinite simple paths in $\gh$, where two paths 
are equivalent if for every finite subgraph $K\subset \gh$,
there is a connected component of $\gh\setminus K$ that intersects
both paths.

Let $\gh$ be a graph.  For a subgraph $H$, let its {\bf (internal) vertex
boundary} $\bdv H$ be the set of vertices of $H$ that are adjacent to some
vertex not in $H$.
We say that $\gh$ is {\bf amenable}
if there is an exhaustion $H_1\subset H_2\subset
\cdots\subset H_n\subset\cdots$ with $\bigcup_n H_n = \gh$ and
$$
\lim_{n \to\infty} {|\bdv H_n| \over |\verts(H_n)|} = 0\,.
$$
Such an exhaustion (or the sequence of its vertex sets) is called a {\bf
F\o{}lner sequence}.
A finitely generated group is {\bf amenable} if its Cayley graph is amenable.
For example, every finitely generated abelian group is amenable.
For more on amenability of graphs and groups, see \BLPSgip.

A {\bf homomorphism} $\varphi: \gh_1 \to \gh_2$ from one
graph $\gh_1=(\verts_1,\edges_1)$ to another $\gh_2=(\verts_2,\edges_2)$
is a pair of maps $\varphi_\verts:\verts_1\to\verts_2$ and
$\varphi_\edges:\edges_1\to\edges_2$
such that $\varphi_\verts$ maps the endpoints of $e$ to the endpoints
of $\varphi_\edges(e)$ for every edge $e \in \edges_1$.
When both maps $\varphi_\verts:\verts_1\to\verts_2$ and
$\varphi_\edges:\edges_1\to\edges_2$ are bijections, then $\varphi$ is
called an {\bf isomorphism}.
When $\gh_1 = \gh_2$, an isomorphism is called
an {\bf automorphism}.
The group of all automorphisms of $\gh$ will be denoted by $\Aut(\gh)$.
The action of a group $\gp$ on a graph $\gh$ by automorphisms is said to be {\bf
transitive} if there is only one $\gp$-orbit in $\verts(\gh)$ and to be {\bf
quasi-transitive} if there are only finitely many orbits in $\verts(\gh)$.
A graph $\gh$ is {\bf transitive} or {\bf quasi-transitive} according as
whether the corresponding action of $\Aut(\gh)$ is.
%if for every $x,y\in\verts$, there is an 
%automorphism of $\gh$ taking $x$ to $y$. 
%In other words, the orbit
%space $\verts/\Aut(\gh)$ of $\verts$ under the group $\Aut(\gh)$ of 
%automorphisms of $\gh$ is a singleton.
For example, every Cayley graph is transitive.
%A graph is {\bf quasi-transitive} if
%the orbit space $\verts/\Aut(\gh)$ is finite.

The action of a group on a set is called {\bf free} if the stabilizer of each
element of the set is just the identity element of the group.
For example, every group acts on itself freely by multiplication.

A locally compact group is called {\bf unimodular} if its
left Haar measure is also right invariant. 
In particular, every discrete countable group is unimodular.
We call a graph $\gh$ {\bf unimodular} if $\Aut(\gh)$ is unimodular, where
$\Aut(\gh)$ is given the weak topology generated by its action on $\gh$.
Every Cayley graph and, as \ref b.SoardiWoess/ and \ref b.Salvatori/ proved,
every quasi-transitive amenable graph is unimodular.
See \BLPSgip\ for more details on unimodular graphs.

We now recall some definitions pertaining to entropy.
For simplicity, we restrict our discussion to groups acting on graphs,
which is enough for our purposes.
Suppose first that $\mu$ is a probability measure on a finite or countable
set, $X$.
The {\bf entropy} of $\mu$ is 
$$
\ent(\mu) := - \sum_{x \in X} \mu(x) \log \mu(x)
\,.
$$
For example, suppose that $\mu$ is a probability measure on the Borel sets of
$2^\edges := \{ 0, 1 \}^\edges$ in the product topology. 
As usual, we identify an element of $2^\edges$ with the subset (or
``configuration") of edges where the value 1 is taken.
If $H$ is a finite subgraph of $\gh$, write $\mu \restrict H$ for the
restriction of $\mu$ to the $\sigma$-field generated by the restrictions to
$\edges(H)$.
That is, $\mu \restrict H$ is defined on the finite set $2^{\edges(H)}$ and
has an entropy as above.
%If $X$ is a random variable with law $\mu$, then we also write $\ent(X)$ for
%$\ent(\mu)$.
%It is well known that if $X$ and $Y$ are independent, then 
%$$
%H\big(\Seq{X, Y}\big) = \ent(X) + \ent(Y)
%\,,
%\label e.ent-ind
%$$
%where $\Seq{X, Y}$ is an ordered pair.

We shall twice use the following well-known lemma (see, e.g., Lemma
6.2 of \ref b.BurPem/ and p.~11 of \ref b.Bollobas/ for the bound on the
binomial coefficient sum).

\procl l.compare
Let $Y$ be a finite set and $m$ be a positive integer.
Write $\alpha := m/|Y|$.
Suppose that $\mu$ is a probability measure on $2^Y \times 2^Y$ that is
supported on the set of pairs $(\omega_1, \omega_2)$ with $|\omega_1 \xor
\omega_2| \le m$.
Let $\mu_1$ and $\mu_2$ be the coordinate marginals of $\mu$.
Then 
$$
|\ent(\mu_1) - \ent(\mu_2)| \le \log \sum_{k=0}^m {|Y| \choose k}
\le |Y| \big(-\alpha \log \alpha - (1 - \alpha) \log (1 - \alpha)\big)
\,.
$$
\endprocl

Suppose that the finite group $\gp$ acts on $X$ and preserves the measure
$\mu$.
Then the {\bf entropy} of the pair $(\mu,\gp)$ is 
$$
\ent(\mu,\gp) := |\gp|^{-1} \ent(\mu)
\,.
$$

Finally, suppose that $\gp$ is a countable amenable finitely generated
subgroup of $\Aut(\gh)$.
Let $\mu$ be a probability measure on $2^\edges$ that is preserved by $\gp$.
Let $\Seq{\gp_n}$ be a F\o{}lner sequence in $\gp$ and $H$ be a finite
subgraph of $\gh$ such that $\gp H =\gh$, provided such an $H$ exists.
Then the {\bf (metric) entropy} of the pair $(\mu,\gp)$, also called the
{\bf $\gp$-entropy} of $\mu$, is 
$$
\ent(\mu,\gp) :=
\lim_{n \to\infty} |\gp_n|^{-1} \ent\big(\mu \restrict (\gp_n H) \big)
\,.
$$
This does not depend on the choice of $H$.
See \ref b.OrnWeiss/ for more details on entropy.

\bsection{Asymptotic Complexity}{s.asym}

Recall that $\cp(\gh)$ denotes the complexity of the graph $\gh$, i.e., the
number of spanning trees of $\gh$.
Let $p_k(x;\gh)$ denote the probability that simple random walk on $\gh$
started at $x$ is back at $x$ after $k$ steps.

We begin with a formula for the complexity of finite graphs.

\procl p.finite
Suppose that $\gh$ is a finite connected graph.
Then 
$$
\log \cp(\gh) 
=
-\log (2|\edges(\gh)|) + \sum_{x \in \verts(\gh)} \log \deg_\gh(x)
- \sum_{k \ge 1} 
{1 \over k} \left(\sum_{x \in \verts(\gh)} p_k(x;\gh) - 1 \right)
\,.
$$
\endprocl

\proof
Let $P$ be the transition matrix for simple random walk on $\gh$ and $I$ be
the identity matrix of the same size.
As shown by \ref b.RungeSachs/, we may rewrite \ref e.MTT/ as
$$
\cp(\gh) = {\prod_{x \in \verts(\gh)} \deg_\gh(x) \over \sum_{x \in
\verts(\gh)} \deg_\gh(x)} \detp (I-P)
%\cp(\gh) = {\prod_{x \in \verts(\gh)} \deg_\gh(x) \over |\vertex(\gh)| }
%\detp [I-P]
$$
[the proof follows from looking at the coefficient of $t$ in
$\det\big(I-P-t I\big) = (\det D_\gh)^{-1} \det\big(\Delta_\gh - t
D_\gh\big)$ and using the Matrix-Tree Theorem in its original form with
cofactors].
Since the sum of the degrees of a graph equals twice the number of its edges,
we obtain
$$
\log \cp(\gh)
= 
-\log (2|\edges(\gh)|) + \sum_{x \in \verts(\gh)} \log \deg_\gh(x)
%-\log |\vertex(\gh)| + \sum_{x \in \verts(\gh)} \log \deg_\gh(x)
+ \log \detp (I-P)
\,.
\label e.mtt'
$$
Let $\Lambda$ be the multiset of eigenvalues of $P$ other than 1 (with
multiplicities).
Since $\Lambda \subset \CO{-1, 1}$, we may rewrite the last term of \ref
e.mtt'/ as 
\begineqalno
\log \detp (I-P)
&=
\sum_{\lambda \in\Lambda} \log (1 - \lambda)
=
- \sum_{\lambda \in\Lambda} \sum_{k \ge 1} \lambda^k/k
\cr&=
- \sum_{k \ge 1} \sum_{\lambda \in\Lambda} \lambda^k/k
=
- \sum_{k \ge 1} {1 \over k} (\tr P^k - 1)
\,,
\endeqalno 
where in the last step, we have used the fact that
the eigenvalue $1$ of $P$ has multiplicity 1
since $\gh$ is connected. 
Since $\tr P^k = \sum_{x \in \verts(\gh)} p_k(x;\gh)$, the desired formula
now follows from this and \ref e.mtt'/.
\Qed

\procl r.mc
There is an extension of \ref p.finite/ to any irreducible Markov chain
with transition matrix $P$. In this case, a spanning tree of the associated
directed graph, with all edges leading towards a single vertex, is often
called a ``spanning arborescence". Let $\cp'(P)$ be
the sum over all spanning trees (with all possible roots) of the product of
$P(e)$ over all $e$ in the tree. The analogue of
\ref p.finite/ states that 
$$
\log \cp'(P)
=
- \sum_{k \ge 1}
{1 \over k} \left(\sum_{x \in \verts(\gh)} p_k(x;\gh) - 1 \right)
\,.
$$
This is reminiscent of a formula that appears in \ref b.LindTun/.
In fact, combining their formula with ours, one gets an expression for the
derivative at $1$ of the reciprocal of the so-called stochastic zeta
function of $P$.
\endprocl

A {\bf rooted graph} $(\gh, x)$ is a graph $\gh$ with a distinguished vertex
$x$ of $\gh$, called the {\bf root}.
A {\bf rooted isomorphism} of rooted graphs is an isomorphism of the
underlying graphs that takes the root of one to the root of the other.
We shall use the following notion of random weak convergence
introduced by \ref b.BS:rdl/ and studied
further by \ref b.AS:obj/ and \ref b.AL:uni/.
Given a positive integer $R$, a finite rooted graph $H$, and a probability
distribution $\rtd$ on rooted graphs, let $p(R, H, \rtd)$ denote the
probability that $H$ is rooted isomorphic to the ball of radius $R$ about
the root of a graph chosen with distribution $\rtd$.
If $(\gh, \rtd)$ is a graph with probability distribution $\rtd$ on its
vertices, then $\rtd$ induces naturally a distribution on rooted graphs,
which we also denote by $\rtd$; namely, the probability of $(\gh, x)$ is
$\rtd(x)$.
For a finite graph $\gh$, let $U(\gh)$ denote the distribution of rooted graphs
obtained by choosing a uniform random vertex of $\gh$ as root of $\gh$.
Suppose that $\Seq{\gh_n}$ is a sequence of finite graphs and that $\rtd$ is a
probability measure on rooted infinite graphs; in most practical cases,
$\rtd$
will be induced by a probability distribution on the vertices of a fixed
infinite graph.
We say the {\bf random weak limit} of $\Seq{\gh_n}$ is $\rtd$ if for any
positive integer $R$ and any finite graph $H$, 
we have $\lim_{n \to\infty} p\big(R, H, U(\gh_n)\big) = p(R, H, \rtd)$.
More generally, if $\gh_n$ are random finite graphs, then
we say the {\bf random weak limit} of $\Seq{\gh_n}$ is $\rtd$ if for any
positive integer $R$, any finite graph $H$, and any $\epsilon > 0$, we have
$\lim_{n \to\infty} \Pleft{\big|p\big(R, H, U(\gh_n)\big) - p(R, H, \rtd)\big| >
\epsilon} = 0$.
%where the left-hand side is a function of the random
%$\gh_n$.
Note that only the component of the root matters for convergence to $\rtd$.
Thus, we may and shall assume that $\rtd$ is concentrated on connected graphs.
If $\rtd$ is induced by a distribution on the vertices of a fixed 
transitive graph $\gh$, then the random weak limit depends only on $\gh$
and not on the root.
In this case, we say that the {\bf random weak limit} of $\Seq{\gh_n}$ is $\gh$. 

Given $R > 0$ and a finite graph $\gh$, let $\nu_R(\gh)$ be the distribution
of the number of edges in the ball of radius $R$ about a random vertex of
$\gh$.
Call a collection of finite graphs $\gh$ {\bf tight} if for each $R$, the
collection of corresponding distributions $\nu_R(\gh)$ is tight.
Note that any tight collection of finite graphs 
has a subsequence that possesses a random weak limit.
%If $\rtd_1$ and $\rtd_2$ are probability measures on $\verts(\gh)$ such that for
%any positive integer $R$ and any finite graph $H$, we have $p(R, H, \rtd_1) =
%p(R, H, \rtd_2)$, then we write $\rtd_1 \equiv \rtd_2$.
%If $\gh$ is a regular graph, write $d_\gh$ for its degree.

Define the {\bf expected degree} of a probability measure $\rtd$ on rooted
graphs to be 
$$
\expdeg (\rtd) := \int \deg_\gh(x) \,d\rtd(\gh, x) 
\,.
$$
When the following integral converges, define the {\bf tree
entropy} of $\rtd$ to be
$$
\asym(\rtd)
:=
\int \Big(\log \deg_\gh(x) - \sum_{k \ge 1} {1 \over k} p_k(x;\gh)\Big)
\,d\rtd(\gh, x) 
\,.
$$
Our terminology is justified by Theorems \briefref t.asym/ and \briefref
t.amen-ent/ below.
The integral converges, for example, when $\rtd$ has finite expected
degree, i.e.,
$
\expdeg (\rtd) < \infty
$,
by virtue of the inequality between the arithmetic and geometric means; see
also \ref c.hbound/ and \ref p.positive/ below.
If $\rtd$ is induced by a fixed transitive graph, $\gh$, of degree $d$, we write 
$$
\asym(\gh) := \asym(\rtd)
=
\log d - \sum_{k \ge 1}
{1 \over k} p_k(o;\gh)
$$
for the tree entropy of $\gh$, where $o$ is any vertex of $\gh$.

Our main theorem in this section is the following.
It also suggests thinking of $\asym(\rtd)$
as a normalized logarithm of the determinant of the Laplacian.
See \ref t.trlog/ for a more direct reason for thinking thus.

\procl t.asym
If $\gh_n$ are finite connected graphs with bounded average degree
whose random weak limit is a probability measure $\rtd$ on infinite rooted
graphs, then 
$$
\lim_{n \to\infty} {1 \over |\verts(\gh_n)|} \log \cp(\gh_n)
=
\asym(\rtd)
\,.
$$
The same limit holds in probability
when $\gh_n$ are random with bounded expected average degree.
\endprocl

\procl r.quotient
If there is a homomorphism that is not an isomorphism from one transitive
graph $\gh_1$ onto another $\gh_2$ of the same degree, then it is clear that
$p_k(o;\gh_1) \le p_k(o; \gh_2)$ for every $k$, with strict inequality for
some $k$. Therefore,
$\asym(\gh_1) > \asym(\gh_2)$. In particular, among all
transitive graphs $\gh$ of degree $d$, the maximum of $\asym(\gh)$ is achieved
uniquely for $\gh$ the regular tree of degree $d$. This maximum value is
calculated explicitly in \ref x.wired-free/.
The uniqueness of the maximum, \ref t.asym/, and tightness
imply that if a sequence of finite regular
graphs $\gh_n$ of degree $d$ does not tend to the $d$-regular tree, $T$, then
$\liminf_n |\verts(\gh_n)|^{-1} \log\cp(\gh_n) < \asym(T)$; this is Theorem
4.5 of \ref b.McKay/.
\endprocl

The only difficulty in deducing \ref t.asym/ from \ref p.finite/ is the
interchange of limit and summation. In order to accomplish this, we shall use
the following lemmas.

\procl l.logQ
Let $P$ be the transition matrix of a Markov chain.
For $\alpha \in [0, 1]$, define the transition matrix $Q := \alpha I + (1 -
\alpha) P$.
For a state $x$, let $p_k(x)$ and $q_k(x)$
denote the return probabilities to $x$ after $k$
steps when the Markov chains start at $x$, where the transition matrices
are $P$ and $Q$, respectively.
We have 
$$
\sum_k q_k(x)/k 
=
- \log (1 - \alpha) + \sum_k p_k(x)/k 
\,.
\label e.logQ
$$
\endprocl

\proof
Write $(\cbuldot, \cbuldot)$ for the ordinary inner product in
$\ell^2(V)$, where $V$ is the state space.
Given $x \in V$ and $z \in (0, 1)$, we have 
\begineqalno
\sum_k q_k(x) z^k/k 
&=
-\big([\log (I-z Q)] \II x, \II x\big)
\cr&=
-\Big(\big[\log \big( (1 - z\alpha) I- z (1 -\alpha) P \big)\big] \II x,
\II x\Big)
\cr&=
-\left(\bigg[\log \Big(I-{z ( 1 -\alpha) \over 1-z\alpha}P\Big) \bigg] \II
x, \II x\right) - \log ({1 - z \alpha}) (\II x, \II x)
\cr&=
- \log (1 - z\alpha) 
+ \sum_k p_k(x) \Big({z ( 1 -\alpha) \over 1-z\alpha}\Big)^k{1 \over k}
\,.
\endeqalno
Letting $z \uparrow 1$, we obtain the desired equation.
\Qed

Next we give a universal bound for the
rate of convergence of $p_k(x;\gh)$ to the stationary probability
$\deg_\gh(x)/2|\edges(\gh)|$ for
finite graphs and to 0 for infinite graphs.
In the case of infinite graphs, such a result is first due to \ref b.CKS/.
Our argument is a modification of that of \ref b.TC:survey/ and seems not
to be written anywhere for the case of finite graphs, although there is some
overlap with the treatment of the special case in Example 2.3.1 of \ref
b.SalCos:stflour/.
For sharper bounds that depend on more information about the graph, see,
e.g., Theorems 3.3.11 and 2.3.1 of \ref b.SalCos:stflour/ for the finite
case and \ref b.BCG/ for the infinite case.
See also \ref b.MP:evolve/.

\procl l.rate
Suppose that $Q$ is a transition matrix of a Markov chain that is
reversible with respect to a positive measure $\pi$. If $\pi$ is finite,
then we assume that $\pi$ is normalized to be a probability measure.
Assume that $c := \inf \big\{ \pi(x) Q(x, y) \st x \ne y \hbox{ and }
Q(x, y) > 0 \big\} > 0$ and that $a := \inf_x Q(x, x) > 0$.
For all states $x$ and all $k \ge 0$, we have
$$
\left|{Q^k(x, x) \over \pi(x)} - 1\right|
\le
\min \Big \{ {1 \over a c \sqrt{k+1}}\,,\  
{1 \over 2 a^2 c^2 (k+1)} \Big \}
\label e.finb
$$
if $\pi$ is finite
and 
$$
{Q^k(x, x) \over \pi(x)} \le {1 \over a c \sqrt{k+1}}
\label e.infinb
$$
if $\pi$ is infinite.
\endprocl

\proof
Write $V$ for the state space and $E$ for the set of pairs $(x, y)$ where
$Q(x, y) > 0$ and $x \ne y$.
Write
$$
c_2(x, y) := \pi(x) Q^2(x, y)
$$
and note that for $(x, y) \in E$, we have 
$$
c_2(x, y) 
\ge \pi(x) [Q(x, x) Q(x, y) + Q(x, y) Q(y, y)] \ge 2 a c
\,.
$$
Define the inner product 
$
\ip{f_1, f_2} := \sum_x f_1(x) f_2(x) \pi(x)
$ on $\ell^2(V, \pi)$.

The case where $\pi$ is infinite is simpler, so we treat that first.
Let $f$ be a function on $V$ with finite support.
Let $x_0$ be a vertex where $|f|$ achieves its maximum.
Then 
$$
\|f\|_\infty = |f(x_0)|
\le
{1\over2} \sum_{(x, y) \in E} |f(x) - f(y)|
\le
{1\over2} \sum_{x, y \in V} c_2(x, y) |f(x) - f(y)| /(2a c)
\,,
\label e.1
$$
where the factor $1/2$ arises from counting each pair $(x, y)$ in each
order.
Applying \ref e.1/ to the function $f^2$, we obtain 
\begineqalno
(2 a c)^2 \|f\|_\infty^4
&\le
\left({1\over2} \sum_{x, y \in V} c_2(x, y) \big|f(x) - f(y)\big|\cdot
    \big|f(x) + f(y)\big| \right)^2
\cr&\le \left({1\over2} \sum_{x, y \in V} c_2(x, y) [f(x) - f(y)]^2\right)
 \left({1\over2} \sum_{x, y \in V} c_2(x, y)[f(x) + f(y)]^2\right)
\cr &= \bigip{ (I-Q^2) f, f } \bigip{ (I+Q^2) f, f }
\endeqalno
by the Cauchy-Schwarz inequality
and some straightforward algebra.
Therefore, if $\ip{f, f} \le 1$, we have 
$
2 (a c)^2 \|f\|_\infty^4
\le \bigip{ (I-Q^2) f, f }
$.
Apply this inequality to the functions 
$Q^l f$ for $l = 0, \ldots, k$ and sum the resulting
inequalities to obtain 
\begineqalno
(k+1) 2 (a c)^2 \|Q^k f\|_\infty ^4 
&\le
2 (a c)^2 \sum_{l=0}^k \|Q^l f\|_\infty ^4 
\le
\sum_{l=0}^k \bigip{ (I-Q^2) Q^l f, Q^l f }
\cr&=
\sum_{l=0}^k \bigip{ (I-Q^2) Q^{2l} f, f }
=
\bigip{ (I-Q^{2k+2}) f, f }
\le
1
\endeqalno
for $\ip{f, f} \le 1$.
This shows that the norm of $Q^k : \ell^2(V, \pi) \to \ell^\infty(V)$ is
bounded by 
$$
\beta_k := [(2k+2) (a c)^2]^{-1/4}
\,.
$$
The same bound holds for $Q^k : \ell^1(V, \pi) \to \ell^2(V, \pi)$ by
duality.
Therefore,
considering $Q^{2k} = Q^k \circ Q^k$, we find that
the norm of $Q^{2k} : \ell^1(V, \pi) \to \ell^\infty(V)$ is at most
$\beta_k^2$,
while the norm of
$Q^{2k+1} : \ell^1(V, \pi) \to \ell^\infty$ is at most $\beta_k
\beta_{k+1}$.
Applying these inequalities to $f := \II x/\pi(x)$ gives \ref e.infinb/.

The case of finite $\pi$ is quite similar.
The essential difference is that we work with $\ell^2_0(V, \pi)$,
the orthogonal complement of the constants in $\ell^2(V, \pi)$.
Note that $\bfo$ is an eigenfunction of $Q$ and that $\ell^2_0(V, \pi)$ is
invariant under $Q$.
We may still conclude \ref e.1/ for all $f$ that have at least one nonnegative
value and at least one nonpositive value, such as
all $f \in \ell^2_0(V, \pi)$.
%The middle part of the argument wll bifurcate, corresponding to
%the two terms on the right-hand side of \ref e.finb/,
%but the last part of the argument wll converge again.

Take $f \in \ell^2_0(V, \pi)$.
Notice that $\sum_{x, y \in V} c_2(x, y) = \sum_{x \in V} \pi(x) = 1$.
Thus, we have from \ref e.1/ that 
$$
(2 a c)^2 \|f\|_\infty^2
\le
{1\over2} \sum_{x, y \in V} c_2(x, y) [f(x) - f(y)]^2
= \bigip{ (I-Q^2) f, f } 
\,.
\label e.sqbd
$$
Alternatively, we may apply \ref e.1/ to the function $\sgn(f)
f^2$.
Using the trivial inequality
$$
|\sgn(s) s^2 - \sgn(t) t^2| \le |s - t| \cdot \big(|s| + |t|\big)
\,,
$$
valid for any real numbers $s$ and $t$, we obtain that
\begineqalno
(2 a c)^2 \|f\|_\infty^4
&\le
\left({1\over2} \sum_{x, y \in V} c_2(x, y) \big|f(x) - f(y)\big|\cdot
    \big(|f(x)| + |f(y)|\big) \right)^2
\cr &\le \bigip{ (I-Q^2) f, f } \bigip{ (I+Q^2) |f|, |f| }
\,.
\endeqalno

Putting both these estimates together, we get
$$
2 (a c)^2 \max \{ 2 \|f\|_\infty^2,\, \|f\|_\infty^4 \}
\le \bigip{ (I-Q^2) f, f }
$$
for $\ip{f, f} = 1$.
As before, this implies that
$$
(k+1) 2 (a c)^2 \max \big\{ 2 \|Q^k f\|_\infty ^2,\, \|Q^k f\|_\infty ^4 \big\}
\le
1
\,,
$$
which shows that the norm of $Q^k : \ell^2_0(V, \pi) \to \ell^\infty(V)$ is
bounded by 
$$
\alpha_k 
:=
\min \big\{ [(2 a c)^2 (k+1)]^{-1/2},\, [(2k+2) (a c)^2]^{-1/4}\big\}
\,.
$$
Let $T :
\ell^2(V, \pi) \to \ell^2_0(V, \pi)$ be 
the orthogonal projection $T f := f - \ip{f, \bfo} \bfo$.
Given what we have shown, we see that
the norm of $Q^k T : \ell^2(V, \pi) \to \ell^\infty(V)$ is 
bounded by $\alpha_k$.
By duality, 
the same bound holds for $T Q^k : \ell^1(V, \pi) \to \ell^2(V, \pi)$.
As before, we deduce that
the norm of $Q^k T Q^k : \ell^1(V, \pi) \to \ell^\infty(V)$ is at most
$\alpha_k^2$ and
the norm of $Q^k T Q^{k+1} : \ell^1(V, \pi) \to \ell^\infty$ is at most
$\alpha_k \alpha_{k+1}$.
Applying these inequalities to $f := \II x/\pi(x)$ gives \ref e.finb/.
\Qed

\procl r.rate
In the infinite case, we do not actually need to assume that $a > 0$.
That is,
suppose that $Q$ is a transition matrix of a Markov chain that is
reversible with respect to a positive measure $\pi$. 
Assume that $c := \inf \big\{ \pi(x) Q(x, y) \st x \ne y \hbox{ and }
Q(x, y) > 0 \big\} > 0$.
If $\pi$ is infinite, then
$$
Q^k(x, x)/\pi(x) \le {4 \over c \sqrt{k+1}}
\,.
%\label e.better
$$
To see this, we have only to establish that
$$
|f(x_0)|
\le
(1/c)\sum_{x, y \in V} c_2(x, y) |f(x) - f(y)|
%\,.
%\label e.2
$$
as a substitute for \ref e.1/.
Let $E_2 := \{ (x, y) \st Q^2(x, y) > 0 \}$.
Because $\pi$ is infinite, for any finite set $S$ of states, there is some
$y$ such that $Q(y, S) Q(y, S^c) > 0$, where $Q(y, A) := \sum_{z \in A}
Q(y, z)$.
Note that $Q(y, S) + Q(y, S^c) = 1$.
Either $Q(y, S)$ or $Q(y, S^c)$ is at least $1/2$ and $\pi(y)$ times
the other is at least $c$, whence 
$$
\pi(y) Q(y, S) Q(y, S^c) \ge c/2
\,.
$$
It follows that 
\begineqalno
\sum_{x \in S} \sum_{z \notin S} c_2(x, z)
&=
\sum_{x \in S} \pi(x) \sum_{z \notin S} Q^2(x, z)
=
\sum_{x \in S} \pi(x) \sum_y \sum_{z \notin S} Q(x, y) Q(y, z)
\cr&=
\sum_y \pi(y) Q(y, S) Q(y, S^c)
\ge
c/2
\,.
\endeqalno
That is, for any cutset of edges $e \in E_2$ that separates $x_0$ from
infinity in the graph $(V, E_2)$, the sum of $c_2(e)$ over the cutset is at
least $c/2$.
Therefore, the max-flow min-cut theorem provides a flow $\theta$ from $x_0$
to infinity of value $c/2$ that is bounded by $c_2$ on each edge in $E_2$.
This yields 
\begineqalno
|f(x_0)| c
&=
\bigg| \sum_{(x, y) \in E_2} \big[f(x) - f(y)\big] \theta(x, y) \bigg|
\le
\sum_{(x, y) \in E_2} |f(x) - f(y)| \cdot |\theta(x, y)|
\cr&\le
\sum_{(x, y) \in E_2} |f(x) - f(y)| c_2(x, y)
\,,
\endeqalno
as desired.
(Recall that each pair $(x, y)$ is counted twice in the sum.)
\endprocl

\procl r.mix
The proofs as given show the very same bounds on the more general
quantities $|Q^k(y, x)/\pi(y) - 1|$ or $Q^k(y, x)/\pi(y)$ for all states
$x, y$ and all $k \ge 0$.
\endprocl

\procl c.hbound
If $\rtd$ is a probability measure on infinite rooted graphs with finite
expected degree, then $\asym(\rtd)$ is finite.
\endprocl

\proof
For a graph $\gh$ with transition matrix $P$, let $Q := (I+P)/2$.
Then $Q$ is the transition matrix of the graph $\gh'$ obtained from $\gh$ by
adding loops to each vertex so as to double its degree.
Write $q_k(x;\gh) := p_k(x;\gh')$.
\ref l.logQ/ tells us that
$$
\sum_k q_k(x;\gh)/k 
=
\log 2 + \sum_k p_k(x;\gh)/k 
\,.
$$
In addition, $Q$ is reversible with respect to the measure $x \mapsto
\deg_\gh(x)$, so that \ref l.rate/ applies with $a \ge 1/2$ and $c \ge 1/2$
(equality holds in both cases when $\gh$ has no loops) to yield
$$
q_k(x;\gh)
\le
{4 \over \sqrt{k+1}} \deg_\gh (x)
\,.
$$
Therefore, 
\begineqalno
\asym(\rtd)
&=
\log 2 + \int \left( \log \deg_\gh(x) - \sum_{k \ge 1} q_k(x;\gh)/k \right)
d\rtd(\gh, x)
\cr&\ge
\log 2 + \int \left( \log \deg_\gh(x) - \sum_{k \ge 1} {4 \over k \sqrt{k+1}}
\deg_\gh (x) \right) d\rtd(\gh, x)
\cr&=
\log 2 + \int \log \deg_\gh(x) \,d\rtd(\gh, x) -
\expdeg(\rtd) \sum_{k \ge 1} {4 \over k \sqrt{k+1}}
\,.
\endeqalno
This gives the corollary by the inequality between the arithmetic and
geometric means.
\Qed

The following lemma is well known.

\procl l.geomconv
Suppose that $Y_n$ are real-valued
random variables that converge in distribution to
$Y$ and that $\sup_n \Ebig{|Y_n|} < \infty$. Then for all continuous
functions $f : \R \to
\R$ such that $\lim_{|x| \to\infty} |f(x)|/|x| = 0$, we have $\lim_{n
\to\infty} \Ebig{f(Y_n)} = \Ebig{f(Y)}$.
\endprocl

\proof
The hypotheses easily imply that $f(Y_n)$ form a uniformly integrable set of
random variables. The continuity of $f$ ensures that $f(Y_n)$ converge in
distribution to $f(Y)$. Together, these imply the conclusion.
\Qed

\proofof t.asym
We claim that
the more general second statement of the theorem follows from the first
statement.
The space of rooted-isomorphism classes of graphs is complete, separable and
metrizable (see \ref b.AS:obj/ for some of the details).
Thus, under the hypothesis of the second statement,
we may assume by 
Skorohod's theorem (see, e.g., Theorem 11.7.2 of \ref b.Dudley/) that
$G_n$ are
defined on a common probability space such that a.s., $G_n$ has a random
weak limit $\rtd$.
Therefore, if the first statement holds, then so does the second.

We now prove the first statement.
Double the degree of each vertex in $\gh_n$ by adding loops to give graphs
$\gh'_n$.
These new graphs have transition matrices $Q_n := (I+P_n)/2$, where $P_n$ is
the transition matrix of $\gh_n$.
Furthermore, the random weak limit of $\gh'_n$ is $\rtd'$, where $\rtd'$ is
obtained from $\rtd$ by doubling the degree of each vertex by adding loops.
By \ref l.logQ/, we have 
$
\asym(\rtd') = \asym(\rtd)
$.
Also, $\cp(\gh'_n) = \cp(\gh_n)$, so it suffices to show that 
$$
\lim_{n \to\infty} {1 \over |\verts(\gh'_n)|} \log \cp(\gh'_n)
=
\asym(\rtd')
\,.
$$
Let $d$ be an upper bound for the average degree of $\gh'_n$, i.e., for all
$n$,
$$
2 |\edges(\gh'_n)| \le d |\verts(\gh'_n)|
\,,
\label e.edgebd
$$
so that $|\verts(\gh'_n)|^{-1} \log \big(2 |\edges(\gh'_n)|\big) \to 0$ as $n
\to\infty$.
Since the degree of a random vertex in $\gh'_n$ converges in distribution to
the degree of the root under $\rtd'$, it follows that 
$$
{1 \over |\verts(\gh'_n)|} \sum_{x \in \verts(\gh'_n)} \log\deg_{\gh'_n}(x)
\to
\int \log\deg_\gh(x) \,d\rtd'(\gh, x)
$$
by \ref l.geomconv/ [use there
$f := \log^+$ and $Y_n$ equals the degree of a uniform vertex in $\gh'_n$].
%By virtue of \ref e.edgebd/, we have
%$|\verts(\gh'_n)|^{-1} \log \big(2 |\edges(\gh'_n)|\big) \to 0$ as $n
%\to\infty$.
Thus, in using \ref p.finite/, we have left to show only that 
$$
\lim_{n \to\infty} 
\sum_{k \ge 1}
|\verts(\gh'_n)|^{-1}{1 \over k}
\left(\sum_{x \in \verts(\gh'_n)} p_k(x;\gh'_n) - 1 \right)
=
\int \sum_{k \ge 1} {1 \over k} p_k(x;\gh) \,d\rtd'(\gh, x)
\,.
$$
By definition and the hypothesis, we have for each $k$ that 
$$
\lim_{n \to\infty} 
|\verts(\gh'_n)|^{-1}\left(\sum_{x \in \verts(\gh'_n)} p_k(x;\gh'_n) - 1 \right)
=
\int p_k(x;\gh) \,d\rtd'(\gh, x)
\,.
$$
\ref l.rate/ applies to $\gh'_n$ with stationary probability measure
$x \mapsto \deg_{\gh'_n}(x)/\big[2|\edges(\gh'_n)|\big]$
and constants $a \ge 1/2$, $c \ge 1/[4|\edges(\gh'_n)|]$ to yield
$$
|\verts(\gh'_n)|^{-1}
\left|\sum_{x \in \verts(\gh_n)} p_k(x;\gh'_n) - 1\right|
\le
{4 d \over \sqrt{k + 1}}
\,.
$$
Hence Weierstrass' M-test justifies the interchange of limit and summation
and we are done.
\Qed

\procl r.weight
A similar result holds for weighted graphs. That is, given a graph $\gh$ whose
edges are assigned positive weights, write $\cp(\gh)$ for the sum of the
weights of its spanning trees, where the weight of a spanning tree is the
product of the weights of its edges. Let the weight of a vertex be the sum
of the weights of its incident edges. The random walk corresponding to a
weighted graph has transition probability from $x$ to $y$ equal to the sum of
the weights of the edges joining $x$ to $y$ divided by the weight of $x$.
If a sequence of weighted connected finite graphs
with weights bounded above and away from 0 and with bounded average vertex
weight has a random weak limit $\rtd$ on weighted rooted infinite graphs,
then the conclusion of \ref t.asym/ holds, where $\asym(\rtd)$ is defined
using the weight of the root in place of its degree and by using
the weighted random walk on the limit graph.
\endprocl

We now illustrate some of the consequences of \ref t.asym/.
Our first result concerns the stability of the asymptotic complexity,
for which we prepare with a lemma related to tightness.

Let $B_R(x) = \big(\vertex_R(x), \edges_R(x)\big)$ denote the ball of
radius $R$ about a vertex $x$.

\procl l.neglig
Let $\Seq{\gh_n}$ be a sequence of finite graphs with vertex
subsets $W_n \subset \vertex(\gh_n)$ satisfying $\lim_{n \to\infty}
|W_n|/v_n = 0$, where
$v_n := |\vertex(\gh_n)|$.
For $R > 0$, let 
$$
s_n(R, t) := |\{x \in \vertex(\gh_n) \st
|\vertex_R(x)| > t\}|
\,.
$$
Let also 
$$
w_n(R) := |\{x \in \vertex(\gh_n) \st
\vertex_R(x) \cap W_n \ne \emptyset\}|
\,.
$$
If for each $R > 0$, we have 
$$
\lim_{t \to\infty} \limsup_{n \to\infty} s_n(R, t)/v_n = 0
\,,
\label e.v-tight
$$
then for each $R > 0$, we have
$\lim_{n \to\infty} w_n(R)/v_n = 0$.
\endprocl

\proof
We have that $\vertex_R(x) \cap W_n \ne \emptyset$
iff $x$ lies in a ball of radius $R$ about some vertex of $W_n$.
If we partition $W_n$ in two, one part consisting of
those vertices $x$ with $|\vertex_R(x)| \le t$
and the other part consisting of the rest, then we deduce that
$$
w_n(R) 
\le
t |W_n| + s_n(2R, t)
$$
for any $R > 0$ and $t > 0$. 
Thus,
$$
\limsup_{n \to\infty} w_n(R)/v_n
\le
\limsup_{n \to\infty} s_n(2R, t)/v_n
\,.
$$
If we now let $t \to\infty$, we obtain the desired result.
\Qed

Note that \ref e.v-tight/ holds if $\Seq{\gh_n}$ has a random weak limit.
In fact, \ref e.v-tight/ is just slightly weaker than tightness, since this
condition counts vertices, while tightness counts edges.

\procl c.subgraph
Let $\Seq{\gh_n}$ be a tight sequence of finite connected graphs with bounded
average degree such that
$
\lim_{n \to\infty} |\vertex(\gh_n)|^{-1}\log\cp(\gh_n) = h
.
$
If $\Seq{\gh'_n}$ is a sequence of connected subgraphs of $\Seq{\gh_n}$
such that
$$
\lim_{n \to\infty} |\vertex(\gh_n)|^{-1}
|\{x \in \vertex(\gh'_n) \st \deg_{\gh'_n}(x) = \deg_{\gh_n}(x) \}|
=
1
\,,
$$
then $\lim_{n \to\infty} |\vertex(\gh'_n)|^{-1}\log\cp(\gh'_n) = h$. 
\endprocl

\proof
By taking a subsequence, if necessary, we may assume (by tightness) that
$\Seq{\gh_n}$ has a random weak limit, $\rtd$.
By \ref t.asym/, we have $h = \asym(\rtd)$.
Let $W_n := \{x \in \vertex(\gh'_n) \st \deg_{\gh'_n}(x) \ne
\deg_{\gh_n}(x) \}$.
Then \ref e.v-tight/ holds because of tightness, whence
$\Seq{\gh'_n}$ also has the random weak limit $\rtd$ by \ref l.neglig/.
Hence \ref t.asym/ applies again to give the desired conclusion.
\Qed

We next illustrate the flexibility of \ref t.asym/ by considering
hybrid graphs as follows.

\procl c.hybrid
Let $\Seq{\gh_n}$ and $\Seq{\gh'_n}$ be tight sequences
of finite connected graphs with bounded
average degree such that
$$
\lim_{n \to\infty} |\vertex(\gh_n)|^{-1}\log\cp(\gh_n) = h
\qquad\hbox{and}\qquad
\lim_{n \to\infty} |\vertex(\gh'_n)|^{-1}\log\cp(\gh_n) = h'
\,.
$$
%Suppose that the random weak limit of the bounded average degree finite
%connected graphs
%$\gh_n$ is $\rtd$ and that
%the random weak limit of the bounded average degree finite connected graphs
%$\gh'_n$ is $\rtd'$.
%\msnote{Reformulate to assume that $\log\cp(\gh_n)/|\vertex(\gh_n)| \to h$ and
%$\log\cp(\gh'_n)/|\vertex(\gh'_n)| \to h'$, using tightness to take convergent
%subsequences.}
Suppose that 
$$
\lim_{n \to\infty} {|\verts(\gh_n)|  \over
|\verts(\gh_n)| + |\verts(\gh'_n)|}
= \alpha \in [0, 1]
\,.
$$
Let $H_n$ be formed by connecting disjoint copies of $\gh_n$ and $\gh'_n$ with
$o\big(|\verts(\gh_n)| + |\verts(\gh'_n)|\big)$ edges in any manner that gives
a sequence of connected graphs.
Then 
$$
\lim_{n \to\infty} {1 \over |\verts(H_n)|} \log \cp(H_n)
=
\alpha h + (1-\alpha)h'
%\alpha \asym(\rtd) + (1-\alpha)\asym(\rtd')
\,.
$$
\endprocl

\proof
By taking subsequences, if necessary, we may assume that
$\Seq{\gh_n}$ and $\Seq{\gh'_n}$ have random weak limits, $\rtd$ and $\rtd'$.
By \ref t.asym/, we have $h = \asym(\rtd)$ and $h' = \asym(\rtd')$.
By \ref l.neglig/,
the random weak limit of $H_n$ is $\alpha \rtd + (1-\alpha)\rtd'$.
Thus \ref t.asym/ gives the desired conclusion.
\Qed

Given probability measures $\rtd_n$ and $\rtd$ on rooted graphs,
we say that $\rtd_n$ {\bf converges weakly} to $\rtd$ if
$p(R, H, \rtd_n) \to p(R, H, \rtd)$ as $n \to\infty$ for any
positive integer $R$ and any finite graph $H$. 
It is not hard to show that tree entropy is a continuous functional when
one bounds the expected degree:

\procl p.cnty
If $\rtd_n$ converges weakly to $\rtd$ as $n \to\infty$
with $\sup_n \expdeg (\rtd_n) < \infty$,
then $\asym(\rtd_n) \to \asym(\rtd)$ as $n \to\infty$.
\endprocl

\proof
%This is the same as saying that given $d < \infty$ and $\epsilon > 0$,
%there is some $R < \infty$ and some $\delta > 0$
%such that for all pairs of probability measures
%$\rtd$ and $\rtd'$ on rooted infinite graphs,
%if $\expdeg(\rtd) \le d$, $\expdeg(\rtd') \le d$, and
%$|\sum_H |p(R, H, \rtd) - p(R, H, \rtd')| < \delta$, then
%$|\asym(\rtd) - \asym(\rtd')| < \epsilon$.
%
As in the proof of \ref t.asym/, we may double the degree of each vertex by
adding loops without changing the tree entropies.
By \ref l.geomconv/ and our assumption of bounded expected degree,
we have $\int \log \deg_\gh(x) \,d\rtd_n(\gh, x) \to 
\int \log \deg_\gh(x) \,d\rtd(\gh, x)$.
Weak convergence itself already guarantees that 
$\int p_k(x; \gh) \,d\rtd_n(\gh, x) \to \int p_k(x; \gh) \,d\rtd(\gh, x)$
for each $k$.
The bounded expected degree and
\ref l.rate/ allow us to apply Weierstrass' M-test to get the desired
conclusion.
\Qed

We now give several explicit examples illustrating the use of tree
entropy, beginning with the transitive case.

In order to evaluate the infinite sum appearing in $\asym(\gh)$, the
following integral is sometimes useful.
Let $\grn(z) := \grn(z;\gh) := \sum_{k \ge 0} p_k(o;\gh) z^k$ be the {\bf
return probability generating function} of the graph $\gh$. Then clearly 
$$
\sum_{k \ge 1} {1 \over k} p_k(o; \gh) = \int_0^1 {\grn(z) - 1 \over z} dz
\,.
\label e.intgrn
$$

\procl x.wired-free
For a group $\gp$ with a given generating set, let $\ell(\gp)$ denote the
length of the smallest (nonempty) reduced word in the generating elements that
represents the identity, i.e., the girth of the Cayley graph of $\gp$.
Suppose that $\gp_n$ are finite groups, each generated by $s$ elements, such
that $\lim_{n \to\infty} \ell(\gp_n) =\infty$.
Then the Cayley graphs $\gh_n$ of $\gp_n$ have a random weak limit equal
to the usual Cayley graph $\gh$
of the free group $\gp$ on $s$ letters, i.e., the regular tree of degree $2s$.
By \ref t.asym/, it follows that 
$$
\lim_{n \to\infty} |\gp_n|^{-1} \log \cp(\gh_n) = \asym(\gh)
%\log 2s - \sum_{k \ge 1} {1 \over k} p_k(o;\gh) =: \asym(\gh) 
\,,
$$
independently of the particular choice of $\gp_n$.
To evaluate $\asym(\gh)$, we use the fact that the return series
is
$$
\grn(z;\gh)
=
{1 - s + \sqrt{s^2-(2s-1)z^2} \over 1 - z^2}
\,,
$$
a result of \ref b.Kesten:symm/.
%Therefore,
%$$
%\sum_{k \ge 1} {1 \over k} p_k(o;\gh)
%=
%\int_{z=0}^1 
%\left({1 - s + \sqrt{s^2-(2s-1)z^2} \over 1 - z^2} - 1\right) {dz \over z}
%\,.
%$$
The integrand in \ref e.intgrn/ then has an ``elementary" antiderivative,
which yields 
$$
\asym(\gh) = \log {(2s-1)^{2s-1}  \over [4 s(s-1)]^{s-1} }
\,.
$$
For example, for $s=2$, we find $\asym(\gh) = \log (3/2)^3$.
%Although entropy is not defined for nonamenable group actions, the number
%$\asym(\gh)$ appears to be the analogue of the entropy of $(\wsf_\gh,\gp)$ in light
%of \ref t.amen-ent/ of the next section.
%It would be interesting to have a meaningful calculation that leads 
%more directly from $\wsf_\gh$ to $\asym(\gh)$.
More generally, when $\gh$ is the regular tree of degree $d$, we have 
$$
\grn(z;\gh) =
{2 (d-1) \over d-2+\sqrt{d^2-4(d-1)z^2}}
$$
(see, e.g., Lemma 1.24 of \ref b.Woess:book/),
whence
$$
\asym(\gh) = \log {(d-1)^{d-1}  \over [d(d-2)]^{d/2-1} }
\,.
$$
For example, $\asym(\gh) = \log (4/\sqrt 3)$ if $d=3$.
As we mentioned in the introduction, this calculation of the asymptotic
complexity of regular graphs with girth tending to infinity was first done by
\ref b.McKay/ under additional hypotheses on the graphs. 
For the case of (uniformly) random $d$-regular graphs, where it is easy to see
that they have a random weak limit equal
to the $d$-regular tree, $\gh$, we
obtain that the asymptotic complexity tends in probability to $\asym(\gh)$;
this was also shown by \ref b.McKay/, who showed in \ref b.McKay:random/ that
random regular graphs satisfy his extra hypotheses.
\endprocl

\procl x.giant
The usual Erd\H{o}s-R\'enyi
model of random graphs, $\G(n, p)$, is a graph on $n$ vertices,
each pair of which is connected by an edge with probability $p$,
independently of other edges.
Other language for this is Bernoulli($p$) bond percolation on the complete
graph $K_n$.
Fix $c > 1$.
The well-known fact that the entire graph $\G(n, c/n)$ has a random weak limit
$\PGW(c)$ is proved explicitly in \ref b.Aldous:mc-pgw/, where $\PGW(c)$ is
the law of a rooted Galton-Watson tree with Poisson($c$) offspring
distribution.
It is well known that with probability approaching 1 as $n \to\infty$,
there is a unique connected component, called the
{\bf giant component}, of $\G(n, c/n)$, that has $\Omega(n)$ vertices.
See, e.g., \ref b.Bollobas/.
Also, the giant component has a random weak limit $\PGW^*(c)$, which is
$\PGW(c)$ conditioned on nonextinction.
This limit of the giant component 
is folklore and seems not to be written anywhere. 
Let $f(c) := \asym\big(\PGW^*(c)\big)$.
We also define $\PGW^*(1)$ to be the weak limit of
$\PGW^*(c)$ as $c \downarrow 1$. Since $\PGW^*(1)$ is the
random weak limit of trees (more specifically, of the uniform spanning tree on
$K_n$) by \ref b.Grimmett:ltrees/,
we have $f(1) = 0$ (which also follows from \ref t.zero/ below).
By \ref p.cnty/, $f$ is continuous on $\CO{1,\infty}$.
We wonder whether $f$ is strictly increasing
on $\CO{1,\infty}$ and real
analytic on $(1,\infty)$.
The fact that $f(c) > 0$ for $c > 1$ follows from \ref t.zero/ below,
together with the well-known fact that $\PGW^*(c)$ has infinitely many ends
a.s.\ for $c > 1$.
It would be interesting to see an explicit formula for $f$.
\endprocl

An additional useful tool for calculation is explained in Section 9
of \ref b.Woess:book/. Namely, as explained
there, if we let
$\rad$ be the radius of convergence of $\grn(\cbuldot)$, then there is a
strictly increasing function $\Phi := \Phi_\gh : \big[0, \rad
\grn(\rad)\big) \to \R$ such that 
$$
\grn(z) = \Phi\big(z \grn(z)\big)
\label e.fnleqn
$$
on $\big[0, \rad \grn(\rad)\big)$.
We are grateful to W.~Woess for pointing out to us the remainder of this
paragraph.
In many cases, it is easier to find $\Phi$ than to (solve and) find $\grn$.
Note that from \ref e.fnleqn/, we have 
$$
\grn(1) = \Phi\big(\grn(1)\big)
\,;
\label e.grn1
$$
it can be much easier to solve for $\grn(1)$ and use that in \ref e.Phiform/
below than it is to solve \ref
e.fnleqn/ for $\grn$ and use that in \ref e.intgrn/.
Substitute \ref e.fnleqn/ in the right-hand side of \ref e.intgrn/ to
obtain 
$$
\int_0^1 {\grn(z) - 1 \over z} dz
=
\int_0^1 {\Phi\big(z\grn(z)\big) - 1 \over z} dz
\,.
$$
Now use the change of variable $t := z \grn(z)$. This gives us
\begineqalno
\int_0^1 {\Phi\big(z\grn(z)\big) - 1 \over z} dz
&=
\int_0^{\grn(1)} {\big(\Phi(t) - 1\big) \big(\Phi(t) - t\Phi'(t)\big) \over t
\Phi(t)} dt
\cr&=
\int_0^{\grn(1)} {\Phi(t) - 1 \over t} dt - \Phi\big(\grn(1)\big) + 1 +
\log \Phi\big(\grn(1)\big)
\cr&=
\int_0^{\grn(1)} {\Phi(t) - 1 \over t} dt - \grn(1) + 1 + \log \grn(1)
\,.
\endeqalno
Thus, we have 
$$
\asym(\gh) = \log d -
\int_0^{\grn(1)} {\Phi(t) - 1 \over t} dt + \grn(1) - 1 - \log \grn(1)
\,.
\label e.Phiform
$$

If we are interested only in the asymptotic complexity of finite graphs, for
which $\rtd$ would we want to calculate $\asym(\rtd)$?
This is answered in \ref s.further/, where we shall see that all $\rtd$ whose
underlying graph is a fixed Cayley graph are included, among others.

\procl x.complete
Suppose that $\gh$ is the free product $K_{s_1} * \cdots * K_{s_n}$
of the complete graphs $K_{s_1}$,
\dots, $K_{s_n}$ for some integers $s_j \ge 2$ with $\sum_j s_j \ge 5$.
In other words, $\gh$ is the Cayley graph of the free product of groups of
order $s_j$ with respect to the generating set corresponding to every
element of the factors other than the identities.
In order to calculate $\asym(\gh)$, we shall find it easier to work with the
graph $\gh'$, in which we have added $n$ loops to each vertex of $\gh$.
By \ref l.logQ/, we have $\asym({\gh'}) = \asym(\gh)$.
Let $d := \sum_{j=1}^n s_j$ be the degree of $\gh'$.
Now by \ref b.Woess:free/, we have 
$$
\Phi_{\gh'}(z)
=
1+{z-n \over 2} + {1 \over 2} \sum_{j=1}^n \sqrt{(1-s_j z/d)^2 + 4z/d}
\,.
$$
For example, if $n=2$, then \ref e.grn1/ gives that $\grn(1, \gh') = s_1
s_2 / (s_1 s_2 - s_1 - s_2)$. We then find via \ref e.Phiform/ applied to
$\gh'$ that 
$$
%\asym(\gh) = {1 \over s} \log \left( \Big({2 \over s-2}\Big)^{s-2}
%(s-1)^{2s-2} \right)
\asym(\gh) =
\left(1 - {1 \over s_1}\right) \log(s_1-1) +
\left(1 - {1 \over s_2}\right) \log(s_2-1) +
\left(1 - {1 \over s_1} - {1 \over s_2}\right) \log {s_1+s_2 \over s_1 s_2
- s_1 - s_2}
\,.
$$
For example, if $s_1 = 2$ and $s_2 = 3$, when $\gh$ is a Cayley graph of the
modular group ${\rm PSL}(2, \Z)$ [use the generators 
$\left(\matrix{0&-1\cr 1&0\cr}\right)$ and
$\left(\matrix{0&-1\cr 1&1\cr}\right)$],
%$z \mapsto -1/z$ and $z \mapsto -1/(z+1)$
then $\asym(\gh) = \log\left(2^{2/3}
5^{1/6}\right)$.
As examples of other tree entropies that one may find by similar
means, we mention that
$\asym(\gh) = \log (16/3)$ if $n=3$ with $s_1 = s_2 = s_3 = 3$, while 
$$
\asym(\gh) =
\log(61+9\sqrt{57}) + {1 \over 6} \log(317-33\sqrt{57}) - {7 \over 2} \log
2 - \log 7
= 1.190^+
$$
if $n = 3$ with $s_1 = s_2 = 2$ and $s_3 = 3$.
As an application, suppose that $H_n$ is a random 3-regular 3-uniform
hypergraph on $n$ vertices. 
%NB: there is such a hypergraph for $n \ge 3$ by induction.
Let $G_n$ be the associated graph in which every
hyperedge of $H_n$ is replaced by a clique on its vertices. Then $G_n$ tends
weakly to $K_3 * K_3 * K_3$, whence $n^{-1} \log\cp(G_n) \to \log (16/3)$.
\endprocl

\procl x.free-free
Choose a ray $\Seq{x_m \st m \ge 0}$ in
the regular tree $T$ of degree 3.
Let $\gh_m$ be the ball in $T$ of radius $m$ about $x_m$.
Remove the edge $[x_m, x_{m+1}]$ from $\gh _m$; let
$T_m$ be the connected component of $x_m$ that remains.
Thus, $T_m \subset T_{m+1}$ for all $m$.
Let $\gh := \bigcup_{m \ge 1} T_m$.
The random weak limit of $\Seq{\gh_m}$ is $\rtd$, where $\rtd(\gh, x_m) :=
2^{-m-1}$ for $m \ge 0$.
Since $\cp(\gh_m) = 1$, \ref t.asym/ tells us that $\asym(\rtd) = 0$, i.e.,
$$
\sum_{k \ge 1}
{1 \over k} \sum_{m \ge 0} 2^{-m-1} p_k(x_m; \gh)
= \log \sqrt3
\,.
$$
\endprocl

We comment finally on situations where
the average degree of $\Seq{\gh_n}$ is not bounded.
We suspect that the following holds 
for ``lazy" simple random walk $Q_\gh := (I+P_\gh)/2$
on any {\it simple\/} (unweighted) graph for some universal constant $C$: 
$$
\all {k \ge 2}\qquad
|\vertex(\gh)|^{-1} \big(\tr Q_\gh^k - 1\big)
\le
{C \over \log^2 k}
\,.
\label e.lazysimple
$$
(It may even be true with $\log^2 k$ replaced by something like $k^{1/3}$,
but as Ben Morris has pointed out, nothing better than $k^{1/3}$ is
possible, as shown by the example of two cliques on $n$ vertices joined by
a path of length $n$.)
If \ref e.lazysimple/ holds, this would replace \ref l.rate/
(except in the proof of \ref c.hbound/)
and allow results still more general than \ref t.asym/ for simple graphs.

As we have not been able to establish \ref e.lazysimple/, we consider
instead sequences
$\Seq{\gh_n}$ that are an
expanding family, meaning that the second largest eigenvalue $\lambda_2(\gh_n)$
of $P_{\gh_n}$ is
bounded away from 1. In this case, we do not need \ref l.rate/, since if
$Q_n = Q_{\gh_n}$ is the transition matrix in the proof of \ref t.asym/, we
have 
$$
|\vertex(\gh_n)|^{-1} \big(\tr Q_n^k - 1\big) \le \left({1 +
\lambda_2(\gh_n) \over 2}\right)^k
$$
for all $n, k \ge 1$.
In addition, we do not need 
\ref l.logQ/; that can be replaced with the corresponding statement for
finite Markov chains, 
$$
\sum_{k \ge 1} (\tr Q^k - 1)/k = - \log (1 - \alpha) +
\sum_{k \ge 1} (\tr P^k - 1)/k 
\,,
$$
where $P$ and $Q$ are as in \ref l.logQ/.

If the average degree of $\Seq{\gh_n}$ is unbounded, we must
consider a different normalization of the complexity.
Let us assume that
the limit (as $n \to\infty$) of the return
probability after $k$ steps of simple random walk started at a random
vertex of $\gh_n$ exists; denote it by $p_k$.
For example, if
$\Seq{\gh_n}$ has
a random weak limit $\rtd$ that is concentrated on infinite
graphs of finite degree, then $p_k = \int p_k(\bp; \gh) \,d\rtd(\gh, \bp)$.
For another example,
$p_k = 0$ when $k \ge 1$ for simple graphs whose minimum
degree tends to infinity.
The proof of \ref t.asym/ (as modified above) shows that
$$
\lim_{n \to\infty} 
|\vertex(\gh_n)|^{-1} 
\left[
\log \cp(\gh_n)
- 
\sum_{x \in \vertex(\gh_n)} \log \deg x
\right]
=
- \sum_{k \ge 1} p_k/k 
\,.
$$
For a particular example, if $\gh_n$ is the giant component
of the random graph $\G(n, p_n)$ with $n p_n -
\log n \to \infty$ (which has $n$ vertices with probability tending to 1;
see, e.g., Theorem 9 of Chapter VII in \ref b.Boll:MGT/), one has
$$
\lim_{n \to\infty} 
\big(n^{-1}
\log \cp(\gh_n)
- 
\log (p_n n)\big)
=
0
$$
in probability.
To facilitate comparison to Cayley's theorem (that $\cp(K_n) = n^{n-2}$), we
may state this as $\cp(\gh_n) = \big([1 + o(1)]p_n n\big)^n$.
(The fact that $\Seq{\gh_n}$ is an expanding family a.s.\
is probably folklore. It can be proved as follows:
First, the proof that $\G(n, p_n)$ is connected with probability
approaching 1 is easily modified to show that its isoperimetric constant,
also called conductance, is bounded away from 0 a.s.
Second, a well-known inequality relating this constant to the second
largest eigenvalue gives the result; see, e.g., \ref b.Chung:cheeger/.)
% One first uses Chernoff's bounds for binomial distributions to bound
% $\sum_{k \le n/2} {n \choose k} \P[Bin(k (n-k), p) < a Bin({k \choose 2},
% p)].
% This shows that one can take a number close to $1-p$ in in the proof of
% Bollob\'as by taking $a$ small. There, one gets the exponent $k (1 -
% \omega(n)) - k \log k + (k^2/n) \log n + (k^2/n) \omega(n)$. Since $k =
% o(n)$, we have k^2/n = o(k)$, which gives the bound we want.
%
\comment{
Direct calculation of the eigenvalues similarly shows that
if $\gh_n$ is the binary cube $\{0, 1\}^n$ with nearest-neighbor
edges, then 
$$
\lim_{n \to\infty} 
\big(2^{-n}
\log \cp(\gh_n)
- 
\log n\big)
=
0
\,.
$$
}

%\bsection{Further Study of Tree Entropy}{s.further}
\bsection{Tree Entropy as Log Determinant}{s.further}

Under certain assumptions,
\ref t.asym/ shows that if $\rtd$ is a random weak limit
of finite connected graphs, $\gh_n$, then its tree entropy $\asym(\rtd)$ is a
limit of the logarithm of the determinant of the graph Laplacians of $\gh_n$,
normalized by omitting the zero eigenvalue and by dividing by the number of
vertices of $\gh_n$.
In fact, one may give a formula for $\asym(\rtd)$ directly in terms of a
normalized determinant of the graph Laplacian for infinite graphs.
This is our main purpose in the present section.
We shall also use this formula to prove inequalities for tree entropy and to
calculate easily and quickly the classical tree entropy for Euclidean
lattices.

We first discuss the class of probability measures $\rtd$ to which our formula
will apply.
This class, the class of $\rtd$ that arise as limits of finite graphs, is the
class of unimodular $\rtd$, defined as follows.
Given a rooted graph $(\gh, x)$ and an edge $e$ incident to $x$, define
the involution $\iota(\gh, x, e) := (\gh, y, e)$, where $y$ is the other
endpoint of $e$.
Given a probability measure $\rtd$ on rooted graphs, define the probability
measure $\widehat\rtd$ to be the law of the isomorphism class of
$(\gh, x, e)$, where $(\gh, x)$ is chosen
according to $\rtd$ and $e$ is then chosen uniformly among the edges
incident to $x$.
Also, define $\widetilde \rtd$ to be the (non-probability) measure that is
the result of biasing $\widehat\rtd$ by the degree of the root; that is,
the Radon-Nikod\'ym derivative of $\widetilde \rtd$ with respect to
$\widehat\rtd$ at the isomorphism class of $(\gh, x, e)$ is $\deg_\gh(x)$.
(If the expected degree of the root is finite, one could obtain a
probability measure from $\widetilde \rtd$ by dividing by the expected degree;
but this is not always the case.)
The involution $\iota$ induces a pushforward map $\widetilde\rtd \mapsto
\iota_* \widetilde\rtd$.
We say that $\rtd$ is {\bf unimodular} or {\bf involution invariant} if
$\iota_*\widetilde\rtd = \widetilde\rtd$.
It is easy to see that any $\rtd$ that is a random weak limit of finite
graphs is unimodular, as observed by \ref b.AS:obj/, who introduced the
notion of involution invariance; essentially the same observation occurs in
\ref b.BS:rdl/.
The converse is much harder, but is established in \ref b.AL:uni/.
%The proofs of most of the remaining results in this section
%depend crucially on this converse; we do not know proofs
%that do not.
Intuitively, unimodularity means that, up to isomorphism, all vertices are
equally likely to be the root.
See \ref b.AL:uni/ for a comprehensive treatment of unimodular random
networks.
In particular, it is shown there that a transitive graph is unimodular iff it
is unimodular as a rooted random graph.

The preceding definitions and results extend easily to the class of 
rooted weighted graphs or multi-graphs $(\gh, \bp)$, where
$\gh = \big(\vertex(\gh), \edges(\gh), w\big)$ and $w : \edges(\gh) \to
\CO{0, \infty}$ is a weight function as in \ref r.weight/.
For $x \ne y \in \vertex(\gh)$, let $\Delta_\gh(x, y) := -\sum_e w(e)$,
where the sum is over all the edges between $x$ and $y$, and
$\Delta_\gh(x, x) := \sum_e w(e)$,
where the sum is over all non-loop edges incident to $x$.
We assume that $\Delta_\gh(x, x) < \infty$ for all $x$.
The associated random walk has the transition probability from $x$ to
$y$ of $-\Delta_\gh(x, y)/\Delta_\gh(x, x)$.
Extend the definition of {\bf tree entropy} to probability measures on rooted 
weighted graphs by 
$$
\asym(\rtd)
:=
\int \Big(\log \Delta_\gh(\bp, \bp) - \sum_{k \ge 1} {1 \over k}
p_k(\bp;\gh)\Big)
\,d\rtd(\gh, \bp) 
$$
whenever this integral converges,
where $p_k(\bp; \gh)$ is the return probability after $k$ steps
for the associated random walk.
%in which the transition probability from $x$ to
%$y$ is $-\Delta_\gh(x, y)/\Delta_\gh(x, x)$.

The {\bf (graph) Laplacian} $\Delta_\gh$ just defined determines an operator 
$$
f \mapsto \Big(x \mapsto \sum_{y \in \vertex} \Delta_\gh(x, y) f(y) \Big)
$$
for functions $f$ on $\vertex$ with finite support.
This operator extends by continuity to a bounded linear operator on
all of $\ell^2(\vertex)$ when $\sup_x \Delta_\gh(x, x) < \infty$.
In case we do not have such a uniform bound, we proceed as follows.
Let $U_e$ be a uniform $[0, 1]$-valued random variable
chosen independently for all $e$.
Given $M \in \Z^+$,
let $\gh'$ be the random weighted graph obtained from $\gh$ by letting
the weight of $e$ be $w(e)(1 - U_e/M)$.
Now let $\gh_M$ be the weighted graph formed from $\gh'$
by changing the weight to 0 of 
those edges $e$ whose weights are greater than $M$ or which are not
among the $M$ largest weights of the edges incident to (or equal to) $e$.
Clearly the matrix $\Delta_{\gh_M}$ converges to $\Delta_\gh$ entrywise
a.s.\ as $M \to\infty$.
Since $\Delta_{\gh_M}$ is a bounded self-adjoint positive semi-definite
operator, the operator $\log (\Delta_{\gh_M} +
\epsilon I)$ is bounded for any $\epsilon > 0$, where $I$ denotes
the identity operator on $\ell^2(\vertex)$.
Let $\rtd_M$ be the law of $(\gh_M, \bp)$ when $(\gh, \bp)$ has the law of
$\rtd$.
If $\rtd$ is unimodular, then so is $\rtd_M$.

Now the logarithm of the determinant of a matrix equals the trace of the
logarithm of the matrix.
Furthermore, one usually {\it defines\/} the determinant via this equality
when one has a trace on a von Neumann algebra. 
This is the approach we take.

The trace we use is defined by \ref b.AL:uni/, which we review here.
Suppose that $\rtd$ is a unimodular probability measure on rooted weighted
graphs.
%We omit the weight function $w$ from our notation when convenient.
Let $T : (\gh, \bp) \mapsto T_{\gh, \bp}$ be a measurable assignment of
bounded linear operators $T_{\gh, \bp}: \ell^2\big(\vertex(\gh)\big) \to
\ell^2\big(\vertex(\gh)\big)$ with finite supremum of
the norms $\|T_{\gh, \bp}\|$.
Let $\alg$ be the von Neumann algebra of such operators $T$ that are
equivariant in the sense that for all isomorphisms $\gpe: \gh_1 \to \gh_2$
and all $\bp, x, y \in \vertex(\gh)$, we have
$(T_{\gh_2, \gpe \bp} \II{\gpe x}, \II{\gpe y}) = (T_{\gh_1,
\bp} \II x, \II y)$.
Since $T_{\gh, \bp}$ does not depend on $\bp \in \gh$ for $T \in
\alg$, we shall write $T_\gh$ in place of $T_{\gh, \bp}$ for $T \in \alg$.
%Identify each $x \in \vertex(\gh)$ with the vector $\II x \in
%\ell^2\big(\vertex(\gh)\big)$.
We define the {\bf trace} of $T \in \alg$ to be
$$
\Tr(T) := \Tr_\rtd(T)
:= \Ebig{(T_{\gh} \II \bp, \II \bp)}
:= \int (T_{\gh} \II \bp, \II \bp) \,d\rtd(\gh, \bp)
\,.
$$
For self-adjoint operators $A$ and $B$,
recall that $A \le B$ means that $\big((B-A) v, v) \ge 0$
for all vectors $v$.
Our trace has the following properties:
$\Tr(\cbuldot)$ is linear,
$\Tr(T) \ge 0$ for $T \ge 0$, and $\Tr(S T) = \Tr(T S)$ for $S, T \in
\alg$.
In addition, for any increasing function $f : \R \to \R$ and any
$S \le T$, we have $\Tr \big(f(S)\big) \le \Tr \big(f(T)\big)$.

In effect, we now show the trace formula
$\asym(\rtd) = \log \Det_\rtd\, \Delta_\gh :=
\Tr_\rtd(\log \Delta_\gh)$.
The determinant here is a so-called Fuglede-Kadison determinant; see \ref
b.FugKad/.

\procl t.trlog
If $\rtd$ is a unimodular probability measure on rooted infinite weighted
graphs with 
$$
\int |\log \Delta_\gh(\bp, \bp)| \,d\rtd(\gh, \bp) < \infty
\label e.logfinite
$$
and 
$$
\int {\Delta_\gh(\bp, \bp) \over \inf \{ -\Delta_\gh(x, y) \st x \ne y
\hbox{ and } \Delta_\gh(x, y) \ne 0 \} } d\rtd(\gh, \bp) < \infty
\,,
\label e.bddwts
$$
then 
$$
\asym(\rtd)
=
\lim_{M \to\infty} \lim_{\epsilon \downarrow 0}
\Tr_{\rtd}
\big[\log (\Delta_{\gh_M} + \epsilon I) \big]
%\int \big(\log (\Delta_{\gh_M} + \epsilon I) \bp, \bp\big)
%\,d\rtd(\gh, \bp)
\,.
\label e.trlog
$$
\endprocl

\proof
It is obvious that $\Delta_{\gh_M} + \epsilon I$ is
monotone increasing in $\epsilon$.
Furthermore, it is easy to check that $\Delta_{\gh_M}$ is monotone
increasing in $M$.
Since $\log$ is an increasing function, it follows that
the limits exist in \ref e.trlog/.
%$\big(\log (\Delta_{\gh_M} + \epsilon I) \bp, \bp\big)$ is
%monotone increasing in $\epsilon$.
%Furthermore, it is easy to check that $\Delta_{\gh_M}$ is monotone
%increasing in $M$, whence
%$\big(\log (\Delta_{\gh_M} + \epsilon I) \bp, \bp\big)$ is
%as well.
%This proves that the limits exist in \ref e.trlog/.

The condition \ref e.logfinite/ and Lebesgue's Dominated Convergence
Theorem guarantee that 
$$
\lim_{M \to\infty} \int \log \Delta_\gh(\bp, \bp) \,d\rtd_M(\gh, \bp) 
=
\int \log \Delta_\gh(\bp, \bp) \,d\rtd(\gh, \bp) 
\,.
$$
Clearly, for all $k > 0$, 
$$
\lim_{M \to\infty} \int p_k(\bp; \gh) \,d\rtd_M(\gh, \bp) 
=
\int p_k(\bp; \gh) \,d\rtd(\gh, \bp) 
\,.
$$
By \ref r.rate/ and \ref e.bddwts/, it follows that
$$
\lim_{M \to\infty} \int \sum_{k \ge 1} {1 \over k}
p_k(\bp; \gh) \,d\rtd_M(\gh, \bp) 
=
\int \sum_{k \ge 1} {1 \over k}p_k(\bp; \gh) \,d\rtd(\gh, \bp) 
\,.
$$
Therefore, we have 
$$
\lim_{M \to\infty} \asym(\rtd_M) = \asym(\rtd)
\,.
\label e.Mtinf
$$

Let $D_\gh$ be the diagonal matrix that has the same diagonal as $\Delta_\gh$.
For $\epsilon > 0$, define $P_{\gh, \epsilon}$ by
$$
(D_\gh + \epsilon I)
(I - P_{\gh, \epsilon})
=
\Delta_\gh + \epsilon I
\,.
$$
In other words, $P_{\gh, \epsilon}(x, x) = 0$ and $P_{\gh, \epsilon}(x, y) =
-\Delta_\gh(x, y)/(D_\gh(x, x) + \epsilon)$ for $x \ne y$.
The matrix $P_{\gh, \epsilon}$ defines a killed random walk that, at $x$, is
killed (sent to an absorbing cemetery state outside the graph $\gh$)
with probability $\epsilon/(D_\gh(x, x) + \epsilon)$ and transits to
$y$ with probability $P_{\gh, \epsilon}(x, y)$.
Let $p_k(\bp; \gh, \epsilon)$ be the return probability after $k$ steps
for the killed random walk.
It is clear that as $\epsilon \downarrow 0$, we have $p_k(\bp; \gh, \epsilon)
\uparrow p_k(\bp; \gh)$.
If the degrees of $\gh$ are bounded, then
the norm of $P_{\gh, \epsilon}$ is less than 1, whence
$$
\log (I - P_{\gh, \epsilon})
=
-\sum_{k \ge 1} {1 \over k}P_{\gh, \epsilon}^k
\,;
$$
in particular,
$$
- \sum_{k \ge 1} {1 \over k}p_k(\bp; \gh, \epsilon) 
=
\big(\log (I - P_{\gh, \epsilon}) \II \bp, \II \bp \big) 
\,.
$$
By a theorem of \ref b.FugKad/, we have therefore for any $M$,
\begineqalno
\asym(\rtd_M)
&=
\int (\log D_\gh \II \bp, \II \bp) - \sum_{k \ge 1} {1 \over k}
p_k(\bp; \gh) \,d\rtd_M(\gh, \bp)
\cr&=
\lim_{\epsilon \downarrow 0}
\int \big(\log (D_\gh + \epsilon I) \II \bp, \II \bp\big)
- \sum_{k \ge 1} {1 \over k}p_k(\bp; \gh, \epsilon) \,d\rtd_M(\gh, \bp)
\cr&=
\lim_{\epsilon \downarrow 0}
\int \big(\log (D_\gh + \epsilon I) \II \bp, \II \bp\big)
+ \big(\log (I - P_{\gh, \epsilon}) \II \bp, \II \bp \big) \,d\rtd_M(\gh, \bp)
\cr&=
\lim_{\epsilon \downarrow 0}
\Tr_{\rtd_M}
\big[\log (D_\gh + \epsilon I) 
+ \log (I - P_{\gh, \epsilon}) \big] 
\cr&=
\lim_{\epsilon \downarrow 0}
\Tr_{\rtd_M}
\Big[\log \big((D_\gh + \epsilon I) 
(I - P_{\gh, \epsilon}) \big) \Big]
\cr&=
\lim_{\epsilon \downarrow 0}
\Tr_{\rtd_M}
\big[\log (\Delta_\gh + \epsilon I) \big]
\,.
\endeqalno
Putting together this limit relation with that of \ref e.Mtinf/, we obtain
\ref e.trlog/.
\Qed

\procl r.loops
As $\Delta_\gh$ is unchanged by the addition or deletion of loops, we see
that neither is $\asym(\rtd)$.
\endprocl

As one indication of the usefulness of \ref t.trlog/ beyond its theoretical
interest, we show how it leads immediately to calculation of the
classically known tree entropy for the nearest-neighbor graph on $\Z^d$. 
In this case, the space $\ell^2(\Z^d)$ is isometrically isomorphic to
$L^2\big([0, 1]^d\big)$ (with Lebesgue measure) via the Fourier transform.
Under this isomorphism, the Laplace operator $\Delta_{\Z^d}$ becomes the
operator of multiplication by the function $(s_1, \ldots, s_d) \mapsto 2d -
2\sum_{i=1}^d \cos (2\pi s_i)$, the vector $\II \bp$ becomes the constant
function $\one$, and thus 
$$
\asym(\Z^d) 
=
\int_{[0, 1]^d} \log\left(2d - 2 \sum_{i=1}^d \cos (2\pi s_i)\right) \,ds
\,.
$$

More generally, suppose that $\gh$ is a graph with vertex set $\Z^d \times
K$ for a finite set $K$ and with edge set that is invariant under the
natural action of $\Z^d$.
That is, for each $x \in \Z^d$, there is an $K \times K$ matrix $L^x$ such
that for $x, y \in \Z^d$ and $u, v \in K$, 
$$
\Delta_\gh\big((x, u), (y, v)\big)
=
L^{y-x}(u, v)
\,.
$$
Such graphs $\gh$ are called ``periodic" by \ref b.BurPem/.
Consider the measure $\rtd$ that puts
equal mass on each $\big(\gh, (\bfz, u)\big)$
($u \in K$).
We may regard $\Delta_\gh$ operating on $\ell^2(\Z^d \times K)$ as an
operator $T$ on $\ell^2\big(\Z^d; \ell^2(K)\big)$, that is, on the space of
vector-valued functions
$f : \Z^d \to \ell^2(K)$ with $\sum_{x \in \Z^d} \|f(x)\|^2_{\ell^2(K)} <
\infty$.
Under the Fourier isomorphism with $L^2\big(\T^d; \ell^2(K)\big)$,
the operator $T$ becomes the operator of
multiplication by the matrix-valued function 
$$
M : (s_1, s_2, \ldots, s_d) \mapsto
\sum_{x \in \Z^d} L^x e^{2 \pi i x \cdot s}
\qquad (s = (s_1, s_2, \ldots, s_d) \in [0, 1]^d)
\,,
$$
the vector $\II{(\bfz, u)}$ becomes the constant function $\II u \in
\ell^2(K)$, and thus 
\begineqalno
\asym(\rtd)
&=
|K|^{-1} \sum_{u \in K} \bigpip{(\log \Delta_\gh) \II{(\bfz, u)},
\II{(\bfz, u)}}
\cr&=
|K|^{-1} \sum_{u \in K} \int_{\T^d} \Bigpip{\big(\log M(s)\big) \II u, \II u}
\,ds
\cr&=
|K|^{-1} \int_{\T^d} \tr \big(\log M(s)\big) \,ds 
\cr&=
|K|^{-1} \int_{\T^d} \log \det M(s) \,ds 
\,.
\endeqalno
This is (a slightly simpler version of) the formula in Theorem 6.1(b) of
\ref b.BurPem/.

We consider next some inequalities.
If $\rtd_1$ and $\rtd_2$ are two probability measures on rooted weighted
graphs, let us
say that $\rtd_1$ is {\bf edge dominated} by $\rtd_2$ if there exists a
probability measure $\nu$ on rooted graphs $(\gh, \bp)$ with two sets of
weights $(w_1, w_2)$ such that for all edges $e$, we have $w_1(e)
\le w_2(e)$ and such that the law of $(\gh, \bp, w_i)$ is $\rtd_i$ for
$i = 1, 2$.
We call $\nu$ a {\bf monotone coupling} of $\rtd_1$ and $\rtd_2$.
When the weight of an edge is 0, one can regard it as being absent.

\procl t.domin
If $\rtd_1 \ne \rtd_2$ are unimodular probability measures on rooted weighted
infinite loop-less
graphs that both satisfy \ref e.logfinite/ and \ref e.bddwts/
and $\rtd_1$ is edge dominated
by $\rtd_2$, then $\asym(\rtd_1) < \asym(\rtd_2)$.
\endprocl

The proof relies on the following notion.
A continuous function $f : (0, \infty) \to \R$ is called {\bf operator
monotone on $(0, \infty)$} if for any bounded self-adjoint operators $A, B$
with spectrum in $(0, \infty)$ and $A \le B$, we have $f(A) \le f(B)$.
\ref b.lowner/ proved that
the logarithm is an operator monotone function on $(0, \infty)$
(see also Chapter V of \ref b.Bhatia/).

\proof
As in the proof of \ref t.trlog/, we have that 
$
\Delta_{\gh_M, w_1} 
\le
\Delta_{\gh_M, w_2} 
$.
Since $\log$ is an operator monotone function on $(0, \infty)$, it follows
that 
$\Tr_{\rtd_1}\big(\log(
\Delta_{\gh_M, w_1} + \epsilon I)\big)
\le
\Tr_{\rtd_2}\big(\log(
\Delta_{\gh_M, w_2} + \epsilon I)\big)
$,
so that $\asym(\rtd_1) \le \asym(\rtd_2)$ by \ref t.trlog/.
If $\asym(\rtd_1) = \asym(\rtd_2)$, then by \ref t.trlog/, we have 
$$
\lim_{M \to\infty} \lim_{\epsilon \downarrow 0} 
\left[
\Tr_{\rtd_2}\big(\log(
\Delta_{\gh_M, w_2} + \epsilon I)\big)
-
\Tr_{\rtd_1}\big(\log(
\Delta_{\gh_M, w_1} + \epsilon I)\big)
\right]
=
0
\,.
$$
Let $\nu$ be a monotone coupling of $\rtd_1$ and $\rtd_2$.
Since $\log$ is an operator monotone function on $(0, \infty)$ with $\lim_{t
\downarrow 0} t \log t = 0$, we may apply a result from \ref b.AL:uni/ to
deduce that $\Delta_{\gh_M, w_1} = \Delta_{\gh_M, w_2}$ $\nu$-a.s., i.e., that
$\rtd_1 = \rtd_2$.
\Qed
%As in the proof of \ref t.trlog/, we have that 
%$$
%\Delta_{\gh_M, w_1} + \epsilon \delta_{\gh_M, w_1}
%\le
%\Delta_{\gh_M, w_2} + \epsilon \delta_{\gh_M, w_2}
%\,,
%$$
%whence
%$$
%\log(\Delta_{\gh_M, w_1} + \epsilon \delta_{\gh_M, w_1})
%\le
%\log(\Delta_{\gh_M, w_2} + \epsilon \delta_{\gh_M, w_2})
%\,.
%$$

\procl r.explicit
In case there is a unimodular monotone coupling $\nu$ (via marked graphs) 
of probability measures $\rtd \ne \rtd'$ on rooted infinite
graphs that have finite expected degree, where $\rtd$ is edge dominated
by $\rtd'$, then one can prove an explicit lower bound for the difference
$\asym(\rtd') - \asym(\rtd)$.
As shown in \ref b.AL:uni/, there is then a sequence $\Seq{(\gh_n,\gh'_n)}$
of pairs of finite
connected graphs on the same vertex sets and with the edge set of $\gh_n$
contained in the edge set of $\gh'_n$ such that $\gh_n$ [resp.,
$\gh'_n$] has a random weak limit $\rtd$ [resp., $\rtd'$] with the average
degree of $\gh'_n$ tending to the $\rtd'$-expected degree of the root.
A counting argument then shows
that $\asym(\rtd') - \asym(\rtd) \ge c^2(\log 2)/[d(d+1)^2]$, where $c :=
\nu\big[ (\gh, \gh', \bp) \st \deg_{\gh'}(\bp) \ne \deg_\gh(\bp)\big]$ and
$d := \expdeg(\rtd')$.  %\int \deg_{\gh'}(\bp) \,d\rtd'(\gh', \bp)$.
\endprocl

\ref c.hbound/ shows that $\asym(\rtd) > -\infty$ as long as $\rtd$ has
finite expected degree.
As an example where $\asym(\rtd) < 0$, consider $\rtd$ to be concentrated
on the single rooted graph $(\N, 0)$.
However, this is not possible in the unimodular case 
of unweighted graphs:

\procl p.positive 
If $\rtd$ is a unimodular probability measure on rooted infinite (unweighted)
graphs that
has finite expected degree, then $\asym(\rtd) \ge 0$.
\endprocl

\proof
Under these hypotheses, \ref b.AL:uni/ establish that $\rtd$ is the random
weak limit of a sequence of finite connected graphs with bounded average
degree.
Thus, we may apply \ref t.asym/.
\Qed

Naturally, we wish to know when the tree entropy is $0$.

\comment{
\procl t.domin
If $\rtd \ne \rtd'$ are unimodular probability measures on rooted infinite
graphs that have finite expected degree and $\rtd$ is edge dominated
by $\rtd'$, then $\asym(\rtd) < \asym(\rtd')$.
\endprocl

\proof
As shown in \ref b.AL:uni/, there is a sequence $(\gh_n,\gh'_n)$ of finite
connected graphs on the same vertex sets and with the edge set of $\gh_n$
contained in the edge set of $\gh'_n$ such that $\gh_n$ [resp.,
$\gh'_n$] has a random weak limit $\rtd$ [resp., $\rtd'$] with the average
degree of $\gh'_n$ tending to the $\rtd'$-expected degree of the root.
Let 
$$
A_n := \{ x \in \verts(\gh_n) \st \deg_{\gh'_n}(x) > \deg_{\gh_n}(x) \}
\,.
$$
Since the distributions of the degree of the root are not the same under
$\rtd$ and $\rtd'$,
there is some $c > 0$ such that $|A_n|/|\verts(\gh_n)| > c$
for all large $n$.
Also, there is an upper bound, $d$, on the average degree of $\gh'_n$.
Thus, for each large $n$, there exists a 
set of vertices $x_1, \ldots, x_k \in A_n$
that is independent (i.e., pairwise non-adjacent) in $\gh'_n$
and with $k/| \verts(\gh_n)| > c/(d+1)$.
[For example, assign uniform random labels independently to the vertices in
$A_n$ and choose those whose assigned label is larger than those of all its
neighbors. The chosen set is always independent and has expected proportion
larger than $c/(d+1)$ by Cauchy's inequality, whence there is some choice this
large.]
We now claim that 
$$
\cp(\gh'_n) > 2^{ kc/[d(d+1)]} \cp(\gh_n)
$$
for all large $n$,
so that $\asym(\rtd') \ge c^2\log 2/[d(d+1)^2] + \asym(\rtd)$.
To see this, let $e_i$ be an edge incident to $x_i$ that lies in $\gh'_n$ but
not in $\gh_n$.
For each spanning tree $T$ of $\gh_n$ and each $i$, the addition of $e_i$ to
$T$ creates a cycle through $x_i$, whence there is some edge, say $e'_i$,
incident to $x_i$
in $\gh_n$ that may be deleted from $T \cup \{ e_i \}$ to form a spanning tree
of $\gh'_n$.
Furthermore, because $\{ x_i \}$ is an independent set, for each $i$, we may
simultaneously choose either to keep $e'_i$ or to replace it by $e_i$,
always obtaining a spanning tree of $\gh'_n$, with all these spanning trees
being distinct.
Given two spanning trees $T$ and $T'$ of $\gh_n$, such replacements as we are
considering may lead to the same spanning tree of $\gh'_n$, but this requires
that the same $ \{ e_i \}$ are added.
Such multiple counting leads to the bound 
$$
\cp(\gh'_n) \ge 
\cp(\gh_n) \prod_{i=1}^k \big(1+1/\deg_{\gh_n}(x_i)\big)
\,.
$$
Since $\log (1+t) \ge t \log 2$ for $0 < t \le 1$, we have 
$$
\sum_{i=1}^k \log \Big(1+{1 \over \deg_{\gh_n}(x_i)}\Big)
\ge
\sum_{i=1}^k {\log 2 \over \deg_{\gh_n}(x_i)}
\ge
{k^2 \log 2 \over \sum_{i=1}^k \deg_{\gh_n}(x_i)}
>
{kc \log 2 \over d(d+1)}
\,.
$$
This proves the claim.
\Qed

We may now determine when tree entropy vanishes.
}

\procl t.zero
If $\rtd$ is a unimodular probability measure on rooted infinite (unweighted)
graphs that
has finite expected degree, then $\asym(\rtd) = 0$
iff $\expdeg (\rtd) = 2$ iff $\rtd$-a.s.\ $\gh$ is a
locally finite tree with 1 or 2 ends.
\endprocl

\proof
The last equivalence is proved in \ref b.AL:uni/.
To prove the first equivalence, let $\gh_n$ be finite connected graphs whose
random weak limit is $\rtd$ and with bounded average degree.
Let $T_n$ be a spanning tree of $\gh_n$.
Since $\Seq{\gh_n}$ is a tight sequence, so is $\Seq{T_n}$.
Therefore, by taking a subsequence if necessary, we may assume that
$\Seq{T_n}$ has a random weak limit $\rtd_0$. 
Clearly, $\rtd_0$ is edge dominated by $\rtd$ and
$\asym(\rtd_0) = 0$ (since $\cp(T_n) = 1$).
\ref b.AL:uni/ proved that
$\expdeg (\rtd_0) = 2$.
Thus, $\expdeg (\rtd) = 2$ iff $\rtd = \rtd_0$
iff $\asym(\rtd) = \asym(\rtd_0)$ by \ref t.domin/.
\Qed

\bsection{Metric Entropy}{s.ent}

Suppose that $\gh$ is an infinite quasi-transitive amenable connected graph.
Choose one element $o_i$ from each vertex orbit.
It is shown in \BLPSgip, Proposition 3.6, that there is a probability
measure $\rho$ on the set $ \{ o_i \}$ such that for any F\o{}lner sequence
$\Seq{H_n}$, the relative frequency of vertices in $H_n$ that are in the same
orbit as $o_i$ converges to $\rho(o_i)$. We call this measure $\rho$ the {\bf
natural frequency distribution} of $\gh$.

\procl t.amen-ent
Let $\gh$ be an infinite quasi-transitive amenable connected graph 
with natural frequency distribution $\rho$.
Let $\gh_n$ be finite connected F\o{}lner subgraphs of $\gh$.
Then 
$$
\lim_{n \to\infty} {1 \over |\verts(\gh_n)|} \log \cp(\gh_n)
=
\sum_{x \in \verts(\gh)} \rho(x) \log \deg_\gh(x) - \sum_{k \ge 1}
{1 \over k} \sum_{x \in \verts(\gh)} \rho(x) p_k(x;\gh)
= \asym(\gh, \rho) 
\,.
\label e.top
$$
If $\gp \subseteq \Aut(\gh)$ is a countable finitely generated
group acting freely on
$\verts(\gh)$ with a finite number $I$ of orbits,
then
$$
\ent(\wsf_\gh, \gp) = I \asym(\gh, \rho) \,.
\label e.ent
$$
Furthermore, the $\gp$-entropy of any invariant measure on essential
spanning forests of $\gh$ is at most $I \asym(\gh, \rho)$.
\endprocl

We shall need several lemmas to prove \ref t.amen-ent/. 

\procl l.1end
Let $\gh$ be an infinite quasi-transitive unimodular graph.
If $\gh$ has 2 ends, then $\fo$ is a tree with exactly 2 ends $\wsf$-a.s.,
while otherwise, for $\wsf$-a.e.\ $\fo$, each component tree of $\fo$ has
exactly one end.
\endprocl

\proof
\ref b.Pem:ust/ established this one-endedness for $\gh = \Z^d$
when $2 \le d \le 4$.
\BLPSusf\ proved this in the general case that $\gh$ is transitive (and
unimodular).
This latter proof is long, but not too hard to modify so as to apply to
quasi-transitive graphs.
The case where $\gh$ has 2 ends or is recurrent needs only a few simple
modifications that we do not detail.
The major changes needed in the transient case with other than 2 ends
are as follows.
Let $\{o_i\}$ be a set of representatives of the orbits of $\verts$ and let
$w_i$ be the reciprocal of the Haar measure of the stabilizer of $o_i$,
where we normalize Haar measure so that $\sum_i w_i = 1$.
(In the amenable case, we have $w_i = \rho(o_i)$ by \BLPSgip,
Proposition 3.6.)
Consider any $\Aut(\gh)$-invariant probability measure $\P$ on
$2^{\edges(\gh)}$.
For any subgraph $\omega$ of $\gh$ and any vertex $x$, write $D(x)$ for the
degree of $x$ in $\omega$.
Let $A_i$ be the event that the component of $o_i$ in $\omega$ is
infinite and $p_{\infty, i} := \P[A_i]$.
Then according to Theorem 6.4 of \BLPSgip, we have 
$$
\sum_i w_i \E[D(o_i) \st A_i] \ge \sum_i 2 w_i p_{\infty, i}
\,.
\label e.mindeg
$$
We may use this to prove an analogue of Theorem 7.2 in \BLPSgip, namely,
that if some component of $\omega$ has at least 3 ends with positive
probability, then strict inequality holds in \ref e.mindeg/.
The next step is to combine the proof of Theorem 6.5 of \BLPSusf\ with
Corollary 3.5 of \BLPSgip\ to show that when $\P = \wsf_\gh$,
we have $\sum_i w_i \E[D(o_i)] = 2$, so that equality holds in \ref e.mindeg/.
Therefore, each tree has at most 2 ends $\wsf$-a.s.
The rest of the proof needs simple obvious modifications only, except for the
crucial ``trunk" lemma, i.e., Lemma 10.5 of \BLPSusf.
Almost all of this proof can also be used word for word. The only significant
change needed is that if $o_i$ is on the trunk, then the shift along the trunk
should bring to $o_i$ the next vertex in the orbit of $o_i$ in the direction
of the orientation of the trunk.
\Qed

For our other lemmas, we shall find the following notation convenient.
Given a finite subgraph $H$ of a graph $\gh$ and a configuration $\omega$ of
$\edges(\gh)$, let $H(\omega)$ denote the cylinder event consisting of those
configurations of $\edges(\gh)$ that agree with $\omega$ on $\edges(H)$.
Given also a configuration $\omega$ of $\edges(\gh) \setminus \edges(H)$,
we define two finite graphs from certain vertex identifications on $H$:
Write $H \circ \omega$ for the graph obtained by identifying all vertices of
$H$ that are connected to each other in the graph $\big(\verts(\gh),
\omega\big)$.
Write $H * \omega$ for the graph obtained by identifying all vertices of $H$
that are connected to each other in the graph $\big(\verts(\gh), \omega\big)$
and by identifying all vertices of $H$ that belong to any infinite connected
component in $\big(\verts(\gh), \omega\big)$.
Note that in $H * \omega$, each finite component of $\big(\verts(\gh),
\omega\big)$ yields a separate identification, while the infinite components
of $\big(\verts(\gh), \omega\big)$ yield together one single identification.

The next lemma provides another justification for the adjective ``uniform" in
``wired uniform spanning forest", similar to a Gibbs specification.
However, it does not hold for all graphs.

\procl l.uniform
Let $\gh$ be an infinite quasi-transitive unimodular graph and
let $H$ be a finite connected subgraph of $\gh$.
If $\gh$ has 2 ends, then 
$$
\wsf\Big(H(\fo) \Bigm| \fo \restrict\big(\edges \setminus \edges(H)\big)\Big)
=
\cp\bigg(H \circ
\Big(\fo \restrict\big(\edges \setminus \edges(H)\big)\Big)\bigg)^{-1}
\qquad \wsf\hbox{-a.s.},
$$
while otherwise, 
$$
\wsf\Big(H(\fo) \Bigm| \fo \restrict\big(\edges \setminus \edges(H)\big)\Big)
=
\cp\bigg(H *
\Big(\fo \restrict\big(\edges \setminus \edges(H)\big)\Big)\bigg)^{-1}
\qquad \wsf\hbox{-a.s.}
$$
\endprocl

\proof
The case where $\gh$ has 2 ends is similar, though simpler, than the other
case, so we give the details only for the case where $\gh$ has other than 2
ends.
Let $Z$ be the event that each tree of $\fo$ has exactly one end.
By \ref l.1end/, we have $\wsf(Z) = 1$.

Let $B_R$ denote the ball of radius $R$ about some fixed vertex
of $\gh$.
Choose $R_H$ so that $H \subset B_{R_H}$.
Let $A_R$ be the event that for all $x, y \in \bdv H$ and
$z, w \in \bdv B_R$,
if in $\fo \restrict \big(\edges(B_R) \setminus \edges(H)\big)$

\beginbullets

$x$ is connected to $z$,

$y$ is connected to $w$, and

$x$ is not connected to $w$,

\endbullets

\noindent
then $x$ and $y$ are not connected in $\fo \restrict \edges(H)$.
Thus, $A_R \subseteq A_{R+1}$ for all $R \ge R_H$ and
$$
Z \subseteq \bigcup_R A_R
\,,
$$
whence $\lim_{R \to\infty} \wsf(A_R) = 1$.
%Then for all large $R$, we have 
%$$
%\wsf(A_R) > 1 - \epsilon^2
%\,.
%$$
Since 
$$
\wsf(A_R) = \int \wsf\Big(A_R
\Bigm| \fo \restrict\big(\edges(B_R) \setminus \edges(H)\big)\Big) \,d\wsf
\,,
$$
it follows that for all large $R$, we have 
$$
\wsf(C_R) \ge 1 - \sqrt{1 - \wsf(A_R)}
\,,
%\label e.largeCR
$$
where
$$
C_R := \left\{ \fo \st \wsf\Big(A_R
\Bigm| \fo \restrict\big(\edges(B_R) \setminus \edges(H)\big)\Big) \ge
1 - \sqrt{1 - \wsf(A_R)}\right\}
\,.
$$
In particular, $\wsf(\limsup_{R \to\infty} C_R) = 1$.

By definition,
$$
\wsf\Big(H(\fo) \Bigm| \fo \restrict\big(\edges \setminus \edges(H)\big)\Big)
=
\lim_{R \to\infty}
\wsf\Big(H(\fo) \Bigm| \fo \restrict\big(\edges(B_R)
\setminus \edges(H)\big)\Big)
\quad\wsf\hbox{-a.s.}
$$
Fix a forest $\omega \in Z \cap \limsup_{R \to\infty} C_R$ for which the
limit above holds.
Choose $\epsilon > 0$ arbitrarily small.
Choose $R \ge R_H$ so large that 
$$
\left|
\wsf\Big(H(\omega) \Bigm| \omega \restrict\big(\edges \setminus
\edges(H)\big)\Big)
-
\wsf\Big(H(\omega) \Bigm| \omega \restrict\big(\edges(B_R)
\setminus \edges(H)\big)\Big)\right|
< \epsilon
\,,
\label e.wl
$$
that $\omega \in C_R$, that $\sqrt{1 - \wsf(A_R)} < \epsilon$,
and that each vertex in $\bdv H$ that is connected in $\omega
\restrict \big(\edges(B_R) \setminus \edges(H)\big)$ to $\bdv B_R$ belongs
to an infinite component in $\omega \restrict \big(\edges \setminus
\edges(H)\big)$. 
This last requirement, in combination with $\omega \in Z$, implies that
$\omega \in A_R$.
Consider the cylinder set
$$
D := \big(B_R \setminus \edges(H) \big)(\omega)
=
\Big\{ \fo \st \fo \restrict\big(\edges(B_R) \setminus
\edges(H)\big) = \omega \restrict\big(\edges(B_R) \setminus \edges(H)\big)
\Big \}
\,.
$$
Let $\mu_N$ be the uniform spanning tree measure on $B_N^*$.
By definition,
$$
\wsf\big(H(\omega) \bigm| D \big)
=
\lim_{N \to\infty} \mu_N\big(H(\omega) \bigm| D \big)
$$
and 
$$
\wsf(A_R \mid D)
= \lim_{N \to\infty}
\mu_N(A_R \mid D)
\,.
$$
Since $\omega \in C_R$ and $\sqrt{1 - \wsf(A_R)} < \epsilon$, we have
$\wsf(A_R \mid D) > 1 - \epsilon$.
Thus, we may choose $N > R$ so large that 
$$
\left|\wsf\big(H(\omega) \bigm| D \big)
-
\mu_N\big(H(\omega) \bigm| D \big)\right|
< \epsilon
$$
and $\mu_N(A_R \mid D) > 1 - \epsilon$.
Since 
$$
\mu_N\big(H(\omega) \bigm| D \big)
=
\mu_N(A_R \mid D) \mu_N\big(H(\omega) \bigm| A_R \cap D \big)
+
\mu_N(A_R^c \mid D) \mu_N\big(H(\omega) \bigm| A_R^c \cap D \big)
\,,
$$
we have
$$
\left|\mu_N\big(H(\omega) \bigm| D \big) -
\mu_N\big(H(\omega) \bigm| A_R \cap D \big)\right|
<
2 \epsilon
\,.
$$
Given $A_R \cap D$, the configurations inside $H$ and outside $B_R$ are
$\mu_N$-independent.
Since $\omega \in A_R$, it follows that 
$$
\mu_N\big(H(\omega) \bigm| A_R \cap D \big)
=
\cp(H * \omega)^{-1} 
\,,
$$
and so
$$
\left|\wsf\big(H(\omega) \bigm| D \big) - \cp(H * \omega)^{-1} \right|
<
3 \epsilon
\,.
$$
Therefore,
$$
\left|
\wsf\Big(H(\omega) \Bigm| \omega \restrict\big(\edges \setminus
\edges(H)\big)\Big)
-
\cp(H * \omega)^{-1} \right|
<
4 \epsilon
$$
by \ref e.wl/.
Since $\epsilon$ is arbitrary and $\omega$ is an arbitrary element of a set
of measure 1, the result follows.
\Qed

\procl l.wiring
Let $H$ be a finite connected graph and $W$ be a subset of vertices of $H$.
Let $H'$ be any graph obtained from $H$ by making certain identifications of
the vertices in $W$ with each other.
Write $\alpha := (|W|-1)/|\edges(H)|$.
Then 
$$
|\log \cp(H) - \log \cp(H')| \le |\edges(H)| \big(-\alpha \log \alpha - (1 -
\alpha) \log (1 - \alpha)\big)
\,.
$$
\endprocl

\proof
Let $\mu$ be the uniform spanning tree measure on $H$ and $\mu'$ the uniform
spanning tree measure on $H'$.
It follows from \ref b.FedMih/ that $\mu$ stochastically dominates $\mu'$.
By \refbauthor{Strassen}'s (\refbyear{Strassen}) theorem, this means
that there is a probability measure on pairs $(T, T')$ so that the law
of $T$ is $\mu$, the law of $T'$ is $\mu'$, and $T \supseteq T'$ a.s.
Now $|\edges(T)| = |\verts(H)|-1$ and $|\edges(T')| = |\verts(H')|-1$.
We deduce that a.s.\ 
$$
|\edges(T) \xor \edges(T')| = |\verts(H)| - |\verts(H')| \le |W|-1
\,.
$$
It now follows from \ref l.compare/ that 
$$
|\ent(\mu) - \ent(\mu')| \le |\edges(H)|\big(-\alpha \log \alpha - (1 -
\alpha) \log (1 - \alpha)\big)
\,.
$$
Since $\ent(\mu) = \log \cp(H)$ and $\ent(\mu') = \log \cp(H')$, this is
the same as the desired inequality.
\Qed

\procl l.cylinder
Let $\gh$ be an infinite quasi-transitive unimodular graph and
$H$ be a finite connected subgraph of $\gh$.
Write $\alpha := (|\bdv H|-1)/|\edges(H)|$.
Then for $\wsf$-a.e.\ $\fo$,
$$
|\log \wsf\big(H(\fo)\big)
- \log \cp(H)^{-1}| \le |\edges(H)|\big(-\alpha \log \alpha - (1 -
\alpha) \log (1 - \alpha)\big)
\,.
$$

\endprocl

\proof
Write $\beta := \alpha^{-\alpha} (1 - \alpha)^{-(1-\alpha)}$.
According to \ref l.wiring/, we have $\cp(H')^{-1} \in I_H$, where
$$
 I_H :=
\left[\cp(H)^{-1} \beta^{-|\edges(H)|},\ 
\cp(H)^{-1} \beta^{|\edges(H)|} \right]
$$
and $H' := H \circ \Big(\fo \restrict\big(\edges \setminus
\edges(H)\big)\Big)$ if $\gh$ has 2 ends, while $H' := H * \Big(\fo
\restrict\big(\edges \setminus \edges(H)\big)\Big)$ otherwise.
Combining the preceding relation with \ref l.uniform/, we obtain 
$$
\wsf\big(H(\fo)\big)
=
\Eleft{\wsf\Big(H(\fo) \Bigm|
\fo \restrict\big(\edges \setminus \edges(H)\big)\Big)}
\in
 I_H
$$
a.s. This is the same as the desired inequality.
\Qed

\procl l.essential
Let $H$ be a finite connected graph and $W$ be a subset of vertices of $H$.
Let $\cpe(H, W)$ be the number of spanning forests of $H$ such that each tree
contains at least one vertex of $W$.
Write $\alpha := (|W|-1)/|\edges(H)|$.
Then 
$$
\log \cpe(H, W) \le \log \cp(H) + |\edges(H)|\big(-\alpha \log \alpha - (1 -
\alpha) \log (1 - \alpha)\big)
\,.
$$
\endprocl

\proof
Let $\mu$ be the uniform measure on spanning forests of $H$ such that each
tree contains at least one vertex of $W$.
Let $\nu$ be obtained from $\mu$ by choosing a spanning forest $\fo$ with
distribution $\mu$ and then randomly adding to $\fo$ enough edges of $H$ so
that a spanning tree of $H$ results.
By \ref l.compare/, we have 
$$
|\ent(\mu) - \ent(\nu)| \le |\edges(H)|\big(-\alpha \log \alpha - (1 -
\alpha) \log (1 - \alpha)\big)
$$
since at most $|W|-1$ edges are added.
Since $\ent(\mu) = \log \cpe(H, W)$ and $\ent(\nu) \le \log \cp(H)$, the desired
inequality follows.
\Qed

\proofof t.amen-ent
We first prove that \ref e.top/ holds.
By definition of $\rho$, the graphs $\gh_n$ have a random weak limit $(\gh,
\rho)$.
Thus, \ref e.top/ is a consequence of \ref t.asym/.

We next show \ref e.ent/.
%Write $\wsf := \wsf_\gh$.
Choose a ball $B_R(o)$ of vertices and edges such that $\gp B_R(o) = \gh$.
Let 
$$
\gp_n := \{ \gpe \in \gp \st \gpe B_R(o) \cap \gh_n \ne \emptyset \}
$$
and put 
$$
\gh'_n := \gp_n B_R(o)
\,.
$$
Since $\Seq{\gh_n}$ is a F\o{}lner sequence in $\gh$, it follows that
$\Seq{\gp_n}$ is a F\o{}lner sequence in $\gp$. Therefore,
$$
\ent(\wsf, \gp) =
- \lim_{n \to\infty} |\gp_n|^{-1} \log \wsf\big(\gh'_{n}(\fo)\big)
$$
in $L^1(\wsf)$ by the generalized Shannon-McMillan Theorem of \ref
b.kieffer:sm/.
%for $\wsf$-a.e.\ $\fo$, where $\gh_n := \gp_n B_R(o)$.
Now
$$
\lim_{n \to\infty} |\verts(\gh'_n)|/|\gp_n| = I
$$ 
since $\gp$ acts freely on $\verts$.
Hence
$$
\ent(\wsf, \gp) =
- \lim_{n \to\infty} I |\verts(\gh'_n)|^{-1} \log \wsf\big(\gh'_{n}(\fo)\big)
\label e.sm
$$
in $L^1(\wsf)$.
As we recalled in \ref s.back/, every quasi-transitive amenable graph is
unimodular.
The result now follows from \ref l.cylinder/ in conjunction with \ref
e.sm/.

We finally show that no measure $\mu$ on essential spanning forests of $\gh$
has larger entropy. 
Using some of the same reasoning as above, we have that 
$$
\ent(\mu,\gp) = \lim_{n \to\infty} |\gp_n|^{-1} \ent(\mu \restrict \gh'_n)
\,.
$$
Because $\mu$ is concentrated on essential spanning forests, the number of
elements of the partition generated by $\gh'_n$ that have positive
$\mu$-measure is at most $\cpe(\gh'_n, \bdv\gh'_n)$, whence $\ent(\mu
\restrict \gh'_n) \le \log \cpe(\gh'_n, \bdv\gh'_n)$.
We may apply \ref l.essential/ to obtain the desired
conclusion.
\Qed

\medbreak
\noindent {\bf Acknowledgements.}\enspace I am grateful to Wolfgang Woess,
Thierry Coulhon, Ben Morris,
Yuval Peres, and Scott Sheffield for useful discussions and references. 
Thanks are due to Benny Sudakov for asking about the asymptotics for graphs
whose degree tends to infinity.

\def\noop#1{\relax}
\input \jobname.bbl

\filbreak
\begingroup
\eightpoint\sc
\parindent=0pt\baselineskip=10pt

Department of Mathematics,
Indiana University,
Bloomington, IN 47405-5701
\emailwww{rdlyons@indiana.edu}
{http://mypage.iu.edu/\string~rdlyons/}

\endgroup

\bye